\documentclass[11pt]{article}
\usepackage{amsmath,amsfonts,amssymb,amsthm,bbm,mathptm,cases}
\usepackage{graphicx,graphics,epsfig,subfigure}
\usepackage{graphicx}
\usepackage{booktabs}
\usepackage{threeparttable}
\usepackage{appendix}
\usepackage{mathrsfs}
\usepackage{graphicx}
\usepackage{subfigure}

\usepackage{xcolor,color,soul,array,cite,float,version,times}
\numberwithin{equation}{section}
\numberwithin{figure}{section}
\numberwithin{table}{section}
\numberwithin{footnote}{section}

\theoremstyle{definition}

\newtheorem{sch}{Scheme}[section]
\theoremstyle{plain}
\newtheorem{thm}{Theorem}[section]

\theoremstyle{remark}

\newcommand{\ben}{\begin{eqnarray}}
\newcommand{\een}{\end{eqnarray}}
\newcommand{\bea}{\begin{array}}
\newcommand{\eea}{\end{array}}
\newcommand{\bes}{\begin{subequations}}
\newcommand{\ees}{\end{subequations}}
\newcommand{\bef}{\begin{figure}[H]}
\newcommand{\eef}{\end{figure}}
\newcommand{\bet}{\begin{tikzpicture}}
\newcommand{\eet}{\end{tikzpicture}}
\newcommand{\beq}{\begin{equation}}
\newcommand{\eeq}{\end{equation}}
\def\bena#1\eena{\begin{eqnarray}\begin{array}{l}#1\end{array}\end{eqnarray}}
\def\besl#1\eesl{\begin{subequations}\begin{align}#1\end{align}\end{subequations}}

\newcommand{\parl}[2]{\ensuremath{\frac{\partial #1}{\partial #2}}}

\newcommand{\vparl}[2]{\ensuremath{\frac{\delta #1}{\delta #2}}}

\newcommand{\ohs}[1]{\ensuremath{\overline{#1}^{n+1/2}}}
\newcommand{\hs}[1]{\ensuremath{{#1}^{n+1/2}}}

\def\inc(#1){\includegraphics[height=3 cm]{pics/#1}}

\def\bn{\mathbf{n}}

\def\bx{\mathbf{x}}
\def\br{\mathbf{r}}

\begin{document}
\title{Linear Second Order Energy Stable Schemes of Phase Field  Model with Nonlocal Constraints for Crystal Growth}
\author{
{Xiaobo Jing}\footnote{Beijing Computational Science Research Center, Beijing 100193, P. R. China. } and
{Qi Wang}\footnote{qwang@math.sc.edu, Department of Mathematics,
 University of South Carolina, Columbia, SC 29028, USA;  Beijing Computational Science Research Center, Beijing  100193, P. R. China.}
}
\date{\today}
\maketitle
\begin{abstract}
We present a set of linear, second order, unconditionally energy stable schemes for the Allen-Cahn model with a nonlocal constraint for crystal growth that conserves the mass of each phase. Solvability conditions are established for the linear systems resulting from the linear schemes. Convergence rates are verified numerically.  Dynamics obtained using the   nonlocal Allen-Cahn model are compared with the one obtained  using the classic Allen-Cahn model as well as the Cahn-Hilliard model, demonstrating slower dynamics than that of the  Allen-Cahn model but faster dynamics than that of the Cahn-Hillard model. Thus, the nonlocal Allen-Cahn model can be an alternative   to the Cahn-Hilliard model in simulating crystal growth. Two Benchmark examples are presented to illustrate the prediction made with the nonlocal Allen-Cahn model in comparison to those made with the Allen-Cahn model and the Cahn-Hillard model.

\end{abstract}

\section{Introduction}

\noindent \indent   Phase field crystal (PFC) growth model, developed as an extension to the phase field formalism \cite{elder2002modeling,elder2004modeling,provatas2011phase,jaatinen2010extended, karma1996phase}, has been  successfully applied to various applications in materials science across different time scales \cite{elder2007phase, provatas2007using, asadi2015review}, capturing the interaction between material defects \cite{elder2002modeling} and modeling the microstructure evolution \cite{elder2002modeling,elder2004modeling,elder2007phase, backofen2007nucleation, backofen2010phase, toth2010polymorphism, guttenberg2010emergence, berry2008simulation, lowen2010phase, pisutha2013calculations}. It is a challenge to develop efficient and stable numerical algorithms to faithfully simulate dynamics described by the PFC model. The PFC model is thermodynamically consistent in that the free energy of the thermodynamic model is dissipative. Numerical algorithms that respect the free energy dissipation property at the discrete level are known as energy stable schemes.

 The Cahn-Hilliard equation is a popular phase field model for crystal growth because of its mass (or volume) preserving property. However, the Cahn-Hilliard equation for the crystal growth problem is of up to the 6th order spatial derivative. Searching for a lower order phase field model that can also preserve mass and free energy dissipation properties is therefore a viable alternative.  Allen-Cahn equation is a  popular phase field model which normally has lower spatial derivatives than the Cahn-Hilliard model. It  describes relaxation dynamics of the thermodynamical system to equilibrium. However, in the case of a phase field description,  when the phase variable represents the mass fraction or the volume fraction of a material component, this model does not warrant the conservation of mass of that component. In order to conserve mass, the free energy functional has to be augmented by a mass preserving penalty term or with a   Lagrange multiplier \cite{mellenthin2008phase, yang2006numerical, Rubinstein1992Nonlocal,li2018unconditionally, du2004phase}. This thus modifies the Allen-Cahn equation into a nonlocal equation. We call this the nonlocal Allen-Cahn model or the Allen-Cahn model with a nonlocal constraint.
 Rubinstein and Sternberg studied the Allen-Cahn model with a mass constraint analytically and compared it with the Cahn-Hilliard model \cite{Rubinstein1992Nonlocal}.  Their result seems to favor using the Allen-Cahn model with a mass constraint in place of the Cahn-Hilliard model when studying interfacial dynamics of incompressible, immiscible multi-component material systems.

For the classical Allen-Cahn equation as well as the Cahn-Hilliard equation, there have been several popular numerical approaches to construct energy stable schemes for the equations, including the convex splitting approach \cite{elliott1993global, fan2017componentwise, eyre1998unconditionally, wise2009energy, shen2012second, wang2010unconditionally, han2015second } , the stabilizing approach \cite{shen2010numerical, du2018stabilized, chen2017uniquely}, the energy quadratization (EQ) approach \cite{zhao2018general, yang2017numerical, Gong2018Fully, Gong2018Linear} and the scalar auxiliary variable approach \cite {shen2018scalar, dong2018family, wang2017efficient}.  Recently, the energy quadratization (EQ) and its reincarnation in the scalar auxiliary variable (SAV) method have been applied to a host of thermodynamic and hydrodynamic models owing to their simplicity, ease of implementation, computational efficiency, linearity, and most importantly their energy stability property \cite{Zhaoetal2018-2, yang2017numerical, Gong2018Fully, Gong2018Linear, Yang2017Linear, Zhao2016Numerical, Zhao2017Numerical, chen2018regularized, yang2018efficient, zhao2016energy, zhao2017novel, zhao2016decoupled, dong2018family, wang2017efficient}. Ones have shown that these strategies are general enough to be useful for developing energy stable numerical approximations to any thermodynamically consistent models, i.e., the models satisfy the second law of thermodynamics or are derived from the Onsager principle \cite{onsager1931reciprocal1, onsager1931reciprocal2, zhao2018general}. The   convex splitting approach, scalar auxiliary variable, energy quadratization approach and other methods have been applied to the Cahn-Hillard model for crystal growth \cite{wang2011energy, wise2009energy, shen2018convergence,zhao2018general,yang2017linearly, hu2009stable, elsey2013simple, gomez2012unconditionally, gomez2011provably, tegze2009advanced, vignal2015energy}.

 In this paper,  we develop a set of linear, second order, unconditionally energy stable schemes using the energy quadratization (EQ)  and scalar auxiliary variable (SAV) approach to solve the nonlocal Allen-Cahn equation numerically. The numerical schemes for the Allen-Cahn and the Cahn-Hilliard model are recalled in the paper simply for comparison purposes. In some of these schemes, both EQ and SAV methods are combined to yield linear, energy stable schemes. We note that when a nonlocal Allen-Cahn model is discretized, it is inevitable to yield an integral which has to be treated with a scalar auxiliary variable. When multiple integrals are identified as SAVs in the free energy functional of the Allen-Cahn model with nonlocal constraints, new solution procedures are developed to solve the subproblems in which elliptic equations can be  solved efficiently. All these schemes are linear and second order accurate in time.  On the other hand,when the EQ strategy is coupled with the discretized integrals, the Sherman-Morrison formula can lead us to an efficient numerical scheme as well. In fact, this can be equivalently dealt with using the SAV method, which will be discussed in the Appendix. The numerical schemes developed in this study for the Allen-Cahn equation with  nonlocal  constraints preserve not only mass but also the energy dissipation rate at the discrete level.

 In the end, we conduct two numerical experiments to assess the performance of the schemes. The results  based on EQ and those based on SAV methods perform equally well in preserving mass and the energy dissipation rate.  In addition, the computational efficiency of the schemes is comparatively studied in one of the benchmark examples as well.  To simplify the presentation, we present the temporal discretization of the models using EQ and SAV approaches in detail. Then, we only briefly discuss the strategy to obtain fully discrete schemes by discretizing the semi-discrete schemes in space later. We refer readers to our early publications in \cite{Gong2018Fully,Gong2018Linear} for details. We show that the linear systems resulting from the schemes are all solvable uniquely if the time step size is suitable so that the solution existence and uniqueness in the full-discrete system is warranted.

 The rest of the paper is organized as follows. In $\S 2$, we present the mathematical models for the classical Allen-Cahn, the Cahn-Hilliard, and the Allen-Cahn model with nonlocal constraints. In $\S 3$, we compare their near equilibrium dynamics. In $\S 4$, we present a set of second order, linear, energy stable numerical schemes for the models. In $\S 5$, we conduct mesh refinement tests on all the schemes and carry out two simulations on crystal growth as well as its grain-boundary effects using the models. Finally, we give the concluding remark in section $\S 6$.

\section{Phase Field  Models for Crystal Growth}

\noindent \indent We consider a phase field model for modeling crystal growth in solids with a focus on resolving the detail of transient dynamics. The free energy of the phase field  model for crystal growth is given by \cite{elder2004modeling,jaatinen2010extended,gomez2012unconditionally}
\bena
F=\int_{\Omega} [\frac{\phi}{2}(-\varepsilon+(\nabla ^2+1)^2)\phi+\frac{\phi^4}{4}]\mathrm{d {\bf r}},\label{fe-1}
\eena
where $\phi$ represents an atomistic density field, which is the deviation of the density from the average density and is a conserved field variable, $\epsilon$ is a parameter related to the temperature, that is, higher $\epsilon$ corresponds to a lower temperature, and $\nabla$ is the gradient operator ($\nabla^2$ denotes the Laplacian). In this study, we   use a more general form of the free energy given by
\bena
F=\int_{\Omega} [\frac{\phi}{2}(\nabla ^4+2a \nabla^2 +\alpha)\phi+\frac{\phi^4}{4}]\mathrm{d {\bf r}}\\
=\int_{\Omega} [\frac{1}{2}\|\nabla ^2\phi\|^2+a \phi \nabla^2\phi+\frac{\alpha}{2}\phi^2+\frac{\phi^4}{4}]\mathrm{d {\bf r}}\\
=\int_{\Omega} [\frac{1}{2}(\|\nabla ^2\phi\|+a\phi)^2+\frac{\alpha-a^2}{2}\phi^2+\frac{\phi^4}{4}]\mathrm{d {\bf r}}\\
=\int_{\Omega} [\frac{1}{2}\|\nabla ^2\phi\|^2-a \|\nabla\phi\|^2+\frac{\alpha}{2}\phi^2+\frac{\phi^4}{4}]\mathrm{d {\bf r}},\label{fe-2}
\eena
where $a=1, \alpha=1-\varepsilon$ recovers equation (\ref{fe-1}) and the boundary conditions
\ben
\bn \cdot \nabla \phi=0, \quad \bn \cdot \nabla \nabla^2 \phi=0
\een
are assumed.

Based on the Onsager linear response theory\cite{onsager1931reciprocal1, onsager1931reciprocal2}, transient dynamics of such a system is customarily governed by a time-dependent partial differential equation given by
\bena
\parl{\phi}{t}= -M \mu,\quad \text{in}\, \Omega, \\
\mu=\vparl{F}{\phi}=(\nabla ^2+a)^2\phi+(\alpha-a^2)\phi+\phi^3,
\eena
subject to appropriate boundary and initial conditions, where $M$ and $\mu=\frac{\delta F}{\delta \phi}$ are the mobility matrix and the chemical potential, respectively.

The time rate of change of the free energy is given by
\bena
\frac{dF}{dt}=-\int_{\Omega} \mu M \mu \mathrm{d {\bf r}}+\int_{\partial \Omega} {\bf n}\cdot [(\frac{\partial f}{\partial \nabla \phi}-\nabla \frac{\partial f}{\partial \nabla^2 \phi})\phi_t+\frac{\partial f}{\partial \nabla^2 \phi}\nabla \phi_t] \mathrm d {\bf s}.\label{Energy-dissipation}
\eena
The following boundary conditions will annihilate the boundary terms in the energy dissipation functional
\ben
\bn \cdot \nabla \phi=0, \quad \bn \cdot (\frac{\partial f}{\partial \nabla \phi}-\nabla \frac{\partial f}{\partial \nabla^2 \phi})=0.\label{bc-2}
\een
In the crystal growth model, $\frac{\partial f}{\partial \nabla \phi}=0$.
So, the boundary conditions reduce to
\ben
\bn \cdot \nabla \phi=0, \quad \bn \cdot \nabla \frac{\partial f}{\partial \nabla^2 \phi}=\bn \cdot \nabla \nabla^2 \phi=0. \label{bc-2}
\een
This set of boundary conditions is consistent with the derivation of the free energy functional.

An alternative set of boundary conditions is given by
\ben
\frac{\partial f}{\partial \nabla^2 \phi}=0, \quad \bn \cdot (\frac{\partial f}{\partial \nabla \phi}-\nabla \frac{\partial f}{\partial \nabla^2 \phi})=0.\label{bc-3}
\een
This is equivalent to
\ben
(\nabla^2 +a)\phi=0, \quad \quad \bn \cdot (\frac{\partial f}{\partial \nabla \phi}-\nabla \frac{\partial f}{\partial \nabla^2 \phi})=0.
\een
This  is  different from the previous one.

We can also assign dissipative boundary conditions to the model as follows
\ben
\bn \cdot \nabla \phi_t=0, \quad \phi_t=-\beta\bn \cdot (\frac{\partial f}{\partial \nabla \phi}-\nabla \frac{\partial f}{\partial \nabla^2 \phi}),\label{bc-4}
\een
where $\beta>0$ is inversely proportional to a relaxation time.
The boundary contribution to the energy dissipation is then given by
\ben
-\int_{\partial \Omega} \beta [\bn \cdot (\frac{\partial f}{\partial \nabla \phi}-\nabla \frac{\partial f}{\partial \nabla^2 \phi})]^2 ds.
\een
The total energy dissipation rate or energy dissipation functional is given by
\ben
\frac{dF}{dt}=-\int_{\Omega} \mu M \mu \mathrm{d {\bf r}}-\int_{\partial \Omega} \beta [{\bf n}\cdot (\frac{\partial f}{\partial \nabla \phi}-\nabla \frac{\partial f}{\partial \nabla^2 \phi})]^2 ds.\label{Energy-dissipation}
\een
If $\beta \to \infty$, we recover \eqref{bc-2}.

Two well-known phase field models are the Allen-Cahn and the Cahn-Hilliard equation, whose mobility is given respectively by
\ben
M=\left \{
\bea{ll}
M_0, & \hbox{Allen-Cahn},\\
-\nabla \cdot M_0 \nabla, & \hbox{Cahn-Hilliard},
\eea\right.
\een
where $M_0$ is a prescribed mobility coefficient matrix, which can be a function of $\phi$. The Allen-Cahn equation does not conserve the total mass
$\int_\Omega \phi d\br$ if $\phi$ is the mass-fraction while the Cahn-Hilliard equation does. However, these two models predict  similar near equilibrium dynamics. On the other hand the Allen-Cahn equation is an equation of lower spatial derivatives, and presumably costs less when solved numerically. Thus, one can impose the mass conservation as a constraint to the Allen-Cahn equation for it to be used to describe dynamics in which the mass is conserved. Next, we will briefly recall several ways to enforce mass conservation to dynamics described by the Allen-Cahn equation.

\subsection{Allen Cahn model}

\noindent \indent The classical Allen-Cahn equation  with the non-flux Neumann boundary conditions is given by
\bena
\parl{\phi}{t}= -M \mu,\quad \text{in}\, \Omega,\ \\ \label{AC}
\parl{\phi}{n}=0, \parl{\nabla^2 \phi}{n}=0, \quad \text{in} \, \partial \Omega, \ \\
\phi|_{t=0}=\phi(0, \bx),
\eena
where $M$ is the mobility coefficient and $\mu$ is the chemical potential
\bena
\mu=\vparl{F}{\phi}=(\nabla ^2+a)^2\phi+(\alpha-a^2)\phi+\phi^3.\\
\eena
The energy dissipation rate of the Allen-Cahn equation is given by
\bena
\frac{d F}{dt}=\int_{\Omega} \vparl{F}{\phi}\phi_t \mathrm{d {\bf r}}
=-\int_{\Omega} \mu (M \mu) \mathrm{d {\bf r}} \leq 0,
\eena
provided nonnegative $M$.
The Allen-Cahn model does  not conserve the mass if the mass is denoted as $\int_\Omega \phi \mathrm d{\bf r}$.

\subsection{Cahn-Hillard model}

\noindent \indent The Cahn-Hilliard equation  with the non-flux Neumann boundary condition is given by
\bena
\parl{\phi}{t}= \nabla \cdot (M \nabla \mu),\quad \text{in}\, \Omega,\ \\ \label{CH}
\parl{\phi}{n}=0, \parl{\nabla^2 \phi}{n}=0, \parl{\mu}{n}=0 \quad \text{in} \, \partial \Omega, \ \\
\phi|_{t=0}=\phi(0,\bx),
\eena
where $M$ is the mobility coefficient and $\mu$ is the chemical potential
\bena
\mu=\vparl{F}{\phi}=(\nabla ^2+a)^2\phi+(\alpha-a^2)\phi+\phi^3.\\
\eena
The energy dissipation rate of the equation is given by
\bena
\frac{d F}{dt}=\int_{\Omega} \vparl{F}{\phi}\phi_t \mathrm{d {\bf r}}
=-\int_{\Omega} \nabla\mu M \nabla\mu \mathrm{d {\bf r}} \leq 0,
\eena
provided nonnegative $M$.
The Cahn-Hillard model  conserves the mass. We next discuss the nonlocal Allen-Cahn equation that conserves the mass.

\subsection {Nonlocal Allen-Cahn models}

\noindent \indent we present two methods to impose mass conservation. One is called the Allen-Cahn model with a penalizing potential and the other is called the Allen-Cahn model with a Lagrange multiplier.
\subsubsection{Allen-Cahn model with a penalizing potential}

\noindent \indent In the Allen-Cahn model with a penalizing potential model, a penalizing term is augmented to the free energy to enforce the mass conservation by the model as follows
\bena
F=\int_{\Omega} [\frac{1}{2}\|\nabla ^2\phi\|^2-a \|\nabla\phi\|^2+\frac{\alpha}{2}\phi^2+\frac{\phi^4}{4}]\mathrm{d {\bf r}}+\frac{\eta}{2}(\int_\Omega \phi(t) \mathrm{d {\bf r}}-{\bf M_0})^2,
\eena
where $\eta$ is penalizing parameter, ${\bf M_0}=\int_\Omega \phi(0) \mathrm{d {\bf r}}$ is the initial mass.

The transport equation for $\phi$ is given by the Allen-Cahn equation
\bena
\parl{\phi}{t}= -M \tilde \mu,\label{PP}\\
\parl{\phi}{n}=0,
 \parl{\nabla^2 \phi}{n}=0,   in \; \partial \Omega, \\
 \phi|_{t=0}=\phi(0,\bx)
 \eena
where $M$ is the mobility coefficient and $\tilde \mu$ is the chemical potential given by
\bena
\tilde \mu=\mu+\sqrt{\eta}\zeta,\quad
\zeta=\sqrt{\eta}(\int_\Omega \phi(t) \mathrm{d {\bf r}}-{\bf M_0}).\\
\eena
We calculate the energy dissipation rate as follows
\bena
\frac{d F}{dt}=\int_{\Omega} \vparl{F}{\phi}\phi_t \mathrm{d {\bf r}}
=\int_{\Omega} \tilde \mu (-M \tilde \mu) \mathrm{d {\bf r}} \leq 0,
\eena
provided $M\geq 0$. The modified Allen-Cahn equation is approximately mass preserving depending on the size of $\eta$ and nonlocal. We next discuss another approaches to obtain a nonlocal Allen-Cahn model.

\subsubsection{Allen-Cahn model with a Lagrange multiplier}

\noindent \indent In this model, the free energy is augmented by a penalty term with a Lagrange multiplier $L$ as follows.
\ben
\tilde F=F-L (\int_\Omega \phi(t) \mathrm{d {\bf r}}-{\bf M_0}).
\een
The transport equation for $\phi$ is given by the Allen-Cahn equation
\bena
\parl{\phi}{t}= -M \tilde \mu,\label{La} \\
\parl{ \phi}{n}=0, \parl{\nabla^2 \phi}{n}=0,   in \; \partial \Omega, \\
\phi|_{t=0}=\phi(0,\bx),
\eena
where $M$ is the mobility coefficient and $\tilde \mu$ is the chemical potential given by
\bena
\tilde \mu=\mu-L,
\quad
L=\frac{1}{\int_{\Omega} M  \mathrm{d {\bf r}}}\int_{\Omega}[ M \tilde \mu]\mathrm{d {\bf r}}.
\eena
We calculate the energy dissipation rate as follows
\bena
\frac{d F}{dt}=\int_{\Omega} \vparl{F}{\phi}\phi_t \mathrm{d {\bf r}}
=\int_{\Omega} \tilde \mu (-M \tilde \mu) \mathrm{d {\bf r}} \leq 0,
\eena
provided $M\geq 0$.

\section{Near equilibrium dyamics of the models}

\noindent \indent We study dynamics of the models near their equilibrium solution $\phi^{ss}$. We consider a small perturbation of the steady state given by  $\delta v(t, \br)$:
\bena
\phi=\phi^{ss}+\delta v. \label{PLS}
\eena

For the Allen-Cahn model, substituting (\ref{PLS}) into (\ref{AC}), we get
\bena
\parl{\delta v}{t}= -M ((\nabla ^2+a)^2\delta v-a^2\delta v+(\alpha+3{(\phi^{ss})}^2)\delta v).
\eena
We seek the solution of the linearized partial differential equation system given by
\ben
\delta v= \sum_{k=0,l=0}^{\infty}a_{kl}cos(kx)cos(ly),\label{ansatz}
\een
 in the domain $\Omega=[-\pi,\pi]^2$. Then, we have
\bena
\dot{a}_{kl}(t)=-M a_{kl}(t)[(k^2+l^2)^2-2a(k^2+l^2)+(\alpha+3(\phi^{ss})^2)].
\eena
Instability can emerge if $(k^2+l^2)^2-2a(k^2+l^2)+(\alpha+3(\phi^{ss})^2)<0$, for some wave numbers $k$, $l$.

For the Allen-Cahn model with a penalizing potential, substituting (\ref{PLS}) into the transport equation, we obtain the linearized system as follows
\bena
\parl{\delta v}{t}= -M [(\nabla ^2+a)^2\delta v-a^2\delta v+(\alpha+3{(\phi^{ss})}^2)\delta v+\eta\int_\Omega \delta v \mathrm{d {\bf r}}].
\eena
Using the ansatz \eqref{ansatz}, we have
\bena
\dot{a}_{kl(t)}=-M a_{kl}(t)[(k^2+l^2)^2-2a(k^2+l^2)+\alpha+3(\phi^{ss})^2+4\pi^2\eta\delta_{k0}\delta_{l0}].
\eena
If $(k^2+l^2)^2-2a(k^2+l^2)+\alpha+3(\phi^{ss})^2+4\pi^2\eta\delta_{k0}\delta_{l0}<0$, instability will occur. In comparison, the Allen-Cahn model with a penalizing potential is more stable than the classical Allen-Cahn model.

For the Allen-Cahn model with a Lagrange multiplier, substituting equation \ref{PLS} into the transport equation, we get the linearized system
\bena
\parl{\delta v}{t}= -M [(\nabla ^2+a)^2\delta v-a^2\delta v+(\alpha+3{(\phi^{ss})}^2)\delta v-\frac{\delta v\int_\Omega g(\phi^{ss})\mathrm{d {\bf r}}}{\int_\Omega M \mathrm{d {\bf r}}}-\frac{\int_\Omega g'(\phi^{ss})\delta v\mathrm{d {\bf r}}}{\int_\Omega M \mathrm{d {\bf r}}}+\\
2\frac{\int_\Omega g(\phi^{ss})\mathrm{d {\bf r}}}{(\int_\Omega M \mathrm{d {\bf r}})^2}\int_\Omega M \delta v\mathrm{d {\bf r}}],
\eena
where $g(\phi)=M((\nabla ^2+a)^2\phi+(\alpha-a^2)\phi+\phi^3)$.
Solving the linear system using ansatz \eqref{ansatz}, we have
\bena
\dot{a}_{kl(t)}=-M a_{kl}(t)[(k^2+l^2)^2-2a(k^2+l^2)+\alpha+3{(\phi^{ss})}^2-\frac{\int_\Omega (\alpha+3{\phi^{ss}}^2)\mathrm{d {\bf r}}}{\int_\Omega M\mathrm{d {\bf r}}}\delta_{k0}\delta_{l0}].
\eena
If $(k^2+l^2)^2-2a(k^2+l^2)+\alpha+3{(\phi^{ss})}^2-\frac{\int_\Omega (\alpha+3{(\phi^{ss})}^2)\mathrm{d {\bf r}}}{\int_\Omega M\mathrm{d {\bf r}}}\delta_{k0}\delta_{l0}<0$, instability may ensure. The contribution of the Lagrange multiplier is to introduce a destabilizing mechanism depending on steady state solution $\phi^{ss}$.

For the Cahn-Hillard model, repeating the above analysis, we have the dynamical equation for the Fourier coefficients:
\bena
\dot{a}_{kl(t)}=-M a_{kl}(t)[(k^2+l^2)^2-2a(k^2+l^2)+(\alpha+3(\phi^{ss})^2)](k^2+l^2).
\eena
The window of instability in the Cahn-Hillard model is identical to that in the Allen-Cahn model. However, the growth rates differ.

The linear stability results dictate  initial transient dynamics of the solution towards or away from the given steady state. We will resort to numerical computations for long time transient behavior of the solution.

\section{Numerical Approximations to the Phase Field Models }

\noindent \indent We design  numerical schemes to solve the above nonlocal phase field  equations  to ensure that the energy dissipation property as well as mass conservation are respected.   We do it by  employing the energy quadratization (EQ) and the scalar auxiliary variable method (SAV) developed recently \cite{Zhao2016Numerical, zhao2016energy, shen2018scalar, yang2018efficient}. Both methods depend on a reformulation of the models into equivalent ones with a quadratic energy. From the latter, we have effective ways to design linear  numerical schemes.  For a full review on EQ methods for thermodynamical models, readers are referred to a recent review article \cite{zhao2018general}. All schemes presented below are firstly given as semi-discretized  ones in time and then the full discretization in space will be discussed.
In fact, we have shown recently that BDF and Runge-Kutta methods can be used to design energy stable schemes for thermodynamical systems up to arbitrarily high order in time \cite{Zhaoetal2018-2}.   For comparison purposes, we also present analogous schemes for the classical Allen-Cahn and the Cahn Hilliard model as well.

\subsection{Temporal discretization}

\subsubsection{Numerical methods for the Allen-Cahn model by EQ methods}

\noindent \indent We reformulate the free energy density by introducing an intermediate variable:
\bena
q=\phi^2.
\eena
Then, the free energy recast into
\bena
F
=\int_{\Omega} [\frac{\phi}{2}(\nabla^4+2a\nabla^2+\alpha)\phi+\frac{q^2}{4}]\mathrm{d {\bf r}}.
\eena

We rewrite (\ref{AC}) as
\bena
\parl{\phi}{t}= -M \mu,\quad
\parl{q}{t}= q' \phi_t,\quad
q'=\parl{q}{\phi},
\eena
where
\bena
\mu=\vparl{F}{\phi}=\nabla ^4\phi+2a\nabla ^2\phi+\alpha\phi+\frac{1}{2}qq'.
\eena
We now discretize it using the linear Crank-Nicolson method in time  to arrive at a second order semi-discrete scheme.
\begin{sch} Given initial conditions $\phi^0,q^0$, we first compute $\phi^1, q^1$ by a first order scheme. Having computed $\phi^{n-1},q^{n-1}$, and $\phi^n,q^n$, we compute $\phi^{n+1},q^{n+1}$ as follows.
\bena
\phi^{n+1}-\phi^n= -\Delta t \overline{M}^{n+1/2} [(\hs{\nabla ^4\phi+2a \nabla ^2\phi+\alpha\phi)} +\frac{1}{2}\hs{q} \ohs{q'}],\label{EQAC1}\\
q^{n+1}-q^n= \ohs{q'}(\phi^{n+1}-\phi^n),\label{EQAC3}
\eena
where
\ben
\overline{(\bullet)}^{n+1/2}=\frac{3}{2}(\bullet)^n-\frac{1}{2}(\bullet)^{n-1}, \quad (\bullet)^{n+1/2}=\frac{1}{2}[(\bullet)^{n+1}+(\bullet)^n].
\een
\end{sch}
The numerical implementation can be done as follows
\bena
\phi^{n+1}=A^{-1}b^n,\\
q^{n+1}=q^n+\ohs{q'}(\phi^{n+1}-\phi^n).
\eena
where $A=I+\Delta t\ohs{M} [\frac{\nabla^4}{2}+a\nabla^2+\frac{\alpha}{2}+\frac{{(\ohs{q'})}^2}{4}], b^n=\phi^n-\Delta t \ohs M[\frac{\nabla^4}{2}\phi^n+a\nabla^2\phi^n+\frac{\alpha}{2}\phi^n+\frac{1}{2}q^n\ohs{q'}-\frac{{(\ohs{q'})}^2}{4}\phi^n]$. So, $\phi^{n+1}$ is solved independent of $q^{n+1}$.

We define the discrete energy as follows
\bena
F^n=\int_{\Omega} [\frac{\phi^n}{2}(\nabla^4+2a\nabla^2+\alpha)\phi^{n}+\frac{({q}^{n})^2}{4}]\mathrm{d {\bf r}}.
\eena
\subsubsection{Numerical methods for the Allen-Cahn model by SAV methods}

\noindent \indent Introducing an intermediate variable$r=\sqrt{\int_\Omega \frac{\phi^4}{4}\mathrm d{\bf r}+C_0}$ as the scalar auxiliary variable, the free energy recast into
\bena
F
=\int_{\Omega} [\frac{\phi}{2}(\nabla^4+2a\nabla^2+\alpha)\phi]\mathrm{d {\bf r}}+r^2-C_0.
\eena
We rewrite (\ref{AC}) as
\bena
\parl{\phi}{t}= -M \mu, \mu=\vparl{F}{\phi}=\nabla ^4\phi+2a\nabla ^2\phi+\alpha\phi+2rg,\\
\parl{r}{t}=\int_ \Omega g \phi_t \mathrm d{\bf r}, g=\vparl{r}{\phi}=\frac{\phi^3}{2\sqrt{\int_\Omega\frac{\phi^4}{4}\mathrm d{\bf r}+C_0}}.\\
\eena
We then discretize it using the linear Crank-Nicolson method in time  to arrive at a second order semi-discrete  scheme.
\begin{sch}
 Given initial conditions $\phi^0,r^0$, we first compute $\phi^1, r^1$ by a first order scheme. Having computed $\phi^{n-1},r^{n-1}$, and $\phi^n,r^n$, we compute $\phi^{n+1},r^{n+1}$ as follows.
\bena
\phi^{n+1}-\phi^n= -\Delta t \overline{M}^{n+1/2} \hs \mu,\\
r^{n+1}-r^n=\int_\Omega \ohs g (\phi^{n+1}-\phi^n)\mathrm d{\bf r},
\eena
where
\bena
\hs \mu=\nabla ^4\hs \phi+2a\nabla ^2\hs \phi+\alpha\hs \phi+2\hs r\ohs g.\\
\eena
\end{sch}
We define the discrete energy as follows
\bena
F^n=\int_{\Omega} [\frac{\phi^n}{2}(\nabla^4+2a\nabla^2+\alpha)\phi^{n}]\mathrm{d {\bf r}}+{(r^n)}^2-C_0.
\eena
The numerical scheme can be recast into
\bena
A\phi^{n+1}+(c,\phi^{n+1})d=b^n,\\
r^{n+1}-r^n=\int_\Omega \ohs g (\phi^{n+1}-\phi^n)\mathrm d{\bf r}.
\eena
where $A=I+\Delta t\ohs{M} [\frac{\nabla^4}{2}+a\nabla^2+\frac{\alpha}{2}], c=\ohs g,d=\Delta t\ohs M \ohs g $ and $ b^n=\phi^n-\Delta t \ohs M[\frac{\nabla^4}{2}\phi^n+a\nabla^2\phi^n+\frac{\alpha}{2}\phi^n+2\ohs g r^n-\ohs g \int_\Omega \ohs g \phi^n \mathrm d{\bf r}]$.
Multiplying the inverse of $A$ firstly and taking the inner product of the equation with c secondly, we have
\bena
(c,\phi^{n+1})+(c,\phi^{n+1})(c,A^{-1}d)=(c,A^{-1}b^n).
\eena
Then, the solution in the scheme is solved in the following steps,
\bena
A[x,y]=[d, b^n],\\
\phi^{n+1}=y-\frac{(c,y)}{1+(c,x)}x,\\
r^{n+1}=r^n+\int_\Omega \ohs g (\phi^{n+1}-\phi^n)\mathrm d{\bf r}.
\eena
\subsubsection{Numerical methods for the Cahn-Hillard model by EQ methods}

\noindent \indent The free energy density is reformulated by introducing an intermediate variable:
\bena
q=\phi^2.
\eena
Then, the free energy recast into
\bena
F
=\int_{\Omega} [\frac{\phi}{2}(\nabla^4+2a\nabla^2+\alpha)\phi+\frac{q^2}{4}]\mathrm{d {\bf r}}.
\eena

We rewrite (\ref{CH}) as
\bena
\parl{\phi}{t}= \nabla \cdot (M \nabla \mu),\\
\parl{q}{t}= q' \phi_t.\\
\eena
where
\bena
\mu=\vparl{F}{\phi}=\nabla ^4\phi+2a\nabla ^2\phi+\alpha\phi+\frac{1}{2}qq',\quad
q'=\parl{q}{\phi}.\\
\eena
We then discretize it using the linear Crank-Nicolson method in time  to arrive at a second order semi-discrete  scheme.

\begin{sch} Given initial conditions $\phi^0,q^0$, we first compute $\phi^1, q^1$ by a first order scheme. Having computed $\phi^{n-1},q^{n-1}$, and $\phi^n,q^n$, we compute $\phi^{n+1},q^{n+1}$ as follows.
\bena
\phi^{n+1}-\phi^n= \Delta t \nabla \cdot (\overline{M}^{n+1/2} \nabla[(\hs{\nabla ^4\phi+2a \nabla ^2\phi+\alpha\phi)} +\frac{1}{2}\hs{q} \ohs{q'}]),\label{EQAC1}\\
q^{n+1}-q^n= \ohs{q'}(\phi^{n+1}-\phi^n).\label{EQAC3}
\eena
\end{sch}
The discrete energy is defined as follows
\bena
F^n=\int_{\Omega} [\frac{\phi^n}{2}(\nabla^4+2a\nabla^2+\alpha)\phi^{n}+\frac{({q}^{n})^2}{4}]\mathrm{d {\bf r}}.
\eena
The numerical implementation can be done as follows
\bena
\phi^{n+1}=A ^{-1}b^n,\\
q^{n+1}=q^n+\ohs{q'}(\phi^{n+1}-\phi^n).
\eena
where $A=I-\Delta t\nabla \cdot (\ohs{M} \nabla[\frac{\nabla^4}{2}+a\nabla^2+\frac{\alpha}{2}+\frac{{(\ohs{q'})}^2}{4}]), b^n=\phi^n+\Delta t \nabla \cdot (\ohs M\nabla[\frac{\nabla^4}{2}\phi^n+a\nabla^2\phi^n+\frac{\alpha}{2}\phi^n+\frac{1}{2}q^n\ohs{q'}-\frac{{(\ohs{q'})}^2}{4}\phi^n])$.

\subsubsection{Numerical methods for the Cahn-Hillard model by SAV methods}

\noindent \indent An intermediate variable $r=\sqrt{\int_\Omega \frac{\phi^4}{4}\mathrm d{\bf r}+C_0}$ are introduced to reformulate the free energy, which is

\bena
F=\int_{\Omega} [\frac{\phi}{2}(\nabla^4+2a\nabla^2+\alpha)\phi]\mathrm{d {\bf r}}+r^2-C_0.
\eena
We rewrite (\ref{CH}) as
\bena
\parl{\phi}{t}= \nabla \cdot (M \nabla\mu), \mu=\vparl{F}{\phi}=\nabla ^4\phi+2a\nabla ^2\phi+\alpha\phi+2rg,\\
\parl{r}{t}=\int_ \Omega g \phi_t \mathrm d{\bf r}, g=\vparl{r}{\phi}=\frac{\phi^3}{2\sqrt{\int_\Omega\frac{\phi^4}{4}\mathrm d{\bf r}+C_0}}.\\
\eena
Linear Crank-Nicolson method is used in time  to arrive at a second order semi-discrete  scheme.
\begin{sch}
 Given initial conditions $\phi^0,r^0$, we first compute $\phi^1, r^1$ by a first order scheme. Having computed $\phi^{n-1},r^{n-1}$, and $\phi^n,r^n$, we compute $\phi^{n+1},r^{n+1}$ as follows.
\bena
\phi^{n+1}-\phi^n= \Delta t \nabla \cdot (\overline{M}^{n+1/2} \nabla\hs \mu),\\
r^{n+1}-r^n=\int_\Omega \ohs g (\phi^{n+1}-\phi^n)\mathrm d{\bf r},
\eena
where
\bena
\hs \mu=\nabla ^4\hs \phi+2a\nabla ^2\hs \phi+\alpha\hs \phi+2\hs r\ohs g.\\
\eena
\end{sch}
The discrete energy is defined as follows
\bena
F^n=\int_{\Omega} [\frac{\phi^n}{2}(\nabla^4+2a\nabla^2+\alpha)\phi^{n}]\mathrm{d {\bf r}}+{(r^n)}^2-C_0.
\eena
The numerical scheme can be rewritten into
\bena
A\phi^{n+1}+(c,\phi^{n+1})d=b^n,\\
r^{n+1}-r^n=\int_\Omega \ohs g (\phi^{n+1}-\phi^n)\mathrm d{\bf r}.
\eena
where $A=I-\Delta t\nabla \cdot (\ohs{M} \nabla[\frac{\nabla^4}{2}+a\nabla^2+\frac{\alpha}{2}]), c=\ohs g,d=\Delta t\ohs M \ohs g $ and $ b^n=\phi^n+\Delta t \nabla \cdot (\ohs M\nabla[\frac{\nabla^4}{2}\phi^n+a\nabla^2\phi^n+\frac{\alpha}{2}\phi^n+2\ohs g r^n-\ohs g \int_\Omega \ohs g \phi^n \mathrm d{\bf r}])$.
Multiplying the first equation by the inverse of $A$ firstly and taking the inner product of the equation with c secondly, we have
\bena
(c,\phi^{n+1})+(c,\phi^{n+1})(c,A^{-1}d)=(c,A^{-1}b^n)
\eena
So we   solve the solution in the scheme in the following steps,
\bena
A[x,y]=[d, b^n],\\
\phi^{n+1}=y-\frac{(c,y)}{1+(c,x)}x,\\
r^{n+1}=r^n+\int_\Omega \ohs g (\phi^{n+1}-\phi^n)\mathrm d{\bf r}.
\eena

\subsubsection{Numerical method for the Allen-Cahn model with a penalizing potential by EQ methods}

\noindent \indent In the     Allen-Cahn model with a penalizing potential, we reformulate the free energy density by introducing two intermediate variables
\bena
q=\phi^2,\quad
\zeta=\sqrt{\eta}(\int_\Omega \phi(t) \mathrm{d {\bf r}}-{\bf M_0}).
\eena
Then, the free energy recast into
\bena
F
=\int_{\Omega} [\frac{\phi}{2}(\nabla^4+2a\nabla^2+\alpha)\phi+\frac{q^2}{4}]\mathrm{d {\bf r}}+\frac{\zeta^2}{2}.
\eena
We rewrite the nonlocal Allen-Cahn equation as follows:
\bena
\parl{\phi}{t}= -M \tilde \mu,\\
\parl{\zeta}{t}= \sqrt{\eta}\int_\Omega \parl\phi{t} \mathrm{d {\bf r}},\\
\parl{q}{t}= q' \phi_t.
\eena
where\bena
\tilde u=u+\sqrt{\eta}\zeta,\quad
\mu=\vparl{F}{\phi}=\nabla ^4\phi+2a\nabla ^2\phi+\alpha\phi+\frac{1}{2}qq',\quad
q'=\parl{q}{\phi}.
\eena
We then discretize it using the linear Crank-Nicolson method in time to arrive at a new scheme as follows.
\begin{sch}
Given initial conditions $\phi^0,q^0$, we first compute $\phi^1, q^1$ by a first order scheme. Having computed $\phi^{n-1},q^{n-1}$, and $\phi^n,q^n$, we compute $\phi^{n+1},q^{n+1}$ as follows.
\bena
\phi^{n+1}-\phi^n= -\Delta t \overline{M}^{n+1/2} \hs {\tilde \mu},\label{EQPFAC}\\
\zeta^{n+1}-\zeta^{n}= \sqrt{\eta}\int_\Omega (\phi^{n+1}-\phi^n) \mathrm{d {\bf r}},\\
q^{n+1}-q^n= \ohs{q'}(\phi^{n+1}-\phi^n),
\eena
where
\bena
\hs {\tilde \mu}=(\hs{\nabla ^4\phi+2a \nabla ^2\phi+\alpha\phi)} +\frac{1}{2}\hs{q} \ohs{q'}+\sqrt{\eta}\hs{\zeta}.\label {model23}
\eena
\end{sch}
The discrete energy is defined as follows
\bena
F^n=\int_{\Omega} [\frac{\phi^n}{2}(\nabla^4+2a\nabla^2+\alpha)\phi^{n}+\frac{({q}^{n})^2}{4}]\mathrm{d {\bf r}}+\frac{(\zeta^n)^2}{2}.
\eena
From the scheme, it follows that
\bena
(I+\Delta t \ohs M [\frac{\nabla^4}{2}+a\nabla^2+\frac{\alpha}{2}+\frac{1}{2}{(\ohs{q'})}^2])\phi^{n+1}+\Delta t \ohs M\frac{\eta}{2}\int_\Omega\phi^{n+1}\mathrm{d {\bf r}}=b^n,\\
b^n=(I-\Delta t \ohs M [\frac{\nabla^4}{2}+a\nabla^2+\frac{\alpha}{2}+q^n\ohs{q'}-\frac{1}{2}{(\ohs{q'})}^2])\phi^{n}-\sqrt\eta\zeta^n+\Delta t \ohs M\frac{\eta}{2}\int_\Omega\phi^{n}\mathrm{d {\bf r}}.
\eena
It can be written into a compact form,
\bena
A\phi^{n+1}+(c,\phi^{n+1})d=b^n,
\eena
where $A=I+\Delta t \ohs M [\frac{\nabla^4}{2}+a\nabla^2+\frac{\alpha}{2}+\frac{1}{2}{(\ohs{q'})}^2], c=1, d=\frac{\Delta t \ohs M}{2}$ and $b^n=(I-\Delta t \ohs M [\frac{\nabla^4}{2}+a\nabla^2+\frac{\alpha}{2}+q^n\ohs{q'}-\frac{1}{2}{(\ohs{q'})}^2])\phi^{n}-\sqrt\eta\zeta^n+\Delta t \ohs M\frac{\eta}{2}\int_\Omega\phi^{n}\mathrm{d {\bf r}}. $ \\
Then, the solve is solved  in the following steps,
\bena
A[x,y]=[d, b^n],\\
(c,\phi^{n+1})=\frac{(c,y)}{1+(c,x)},\\
\phi^{n+1}=-(c,\phi^{n+1})x+y.
\eena

\subsubsection{Numerical method for the Allen-Cahn model with a penalizing potential by SAV methods}

\noindent \indent In the penalizing nonlocal Allen-Cahn model, we reformulate the free energy density by introducing two intermediate variables
\bena
r=\sqrt{\int_\Omega \frac{\phi^4}{4} \mathrm{d {\bf r}}+C_0},\quad
\zeta=\sqrt{\eta}(\int_\Omega \phi(t) \mathrm{d {\bf r}}-{\bf M_0}).
\eena
Then, the free energy recast into
\bena
F
=\int_{\Omega} [\frac{\phi}{2}(\nabla^4+2a\nabla^2+\alpha)\phi]\mathrm{d {\bf r}}+r^2-C_0+\frac{\zeta^2}{2}.
\eena
We rewrite the nonlocal Allen-Cahn equation as follows:
\bena
\parl{\phi}{t}= -M \tilde \mu,\\
\parl{\zeta}{t}= \sqrt{\eta}\int_\Omega \parl\phi{t} \mathrm{d {\bf r}},\\
\parl{r}{t}=\int_\Omega g \parl\phi{t} \mathrm{d {\bf r}},
\eena
where\bena
\tilde \mu=\vparl{F}{\phi}=\nabla ^4\phi+2a\nabla ^2\phi+\alpha\phi+2rg+\sqrt{\eta}\zeta,\quad
g=\vparl{r}{\phi}=\frac{\phi^3}{2\sqrt{\int_\Omega\frac{\phi^4}{4}\mathrm d{\bf r}+C_0}}.
\eena
We then discretize it using the linear Crank-Nicolson method in time to arrive at a new scheme as follows.
\begin{sch} Given initial conditions $\phi^0,r^0$, we first compute $\phi^1, r^1$ by a first order scheme. Having computed $\phi^{n-1},r^{n-1}$, and $\phi^n,r^n$, we compute $\phi^{n+1},r^{n+1}$ as follows.
\bena
\phi^{n+1}-\phi^n= -\Delta t \overline{M}^{n+1/2} \hs {\tilde \mu},\label{EQPFAC}\\
\zeta^{n+1}-\zeta^{n}= \sqrt{\eta}\int_\Omega (\phi^{n+1}-\phi^n) \mathrm{d {\bf r}},\\
r^{n+1}-r^n= \int_\Omega \ohs{g}(\phi^{n+1}-\phi^n)\mathrm {d {\bf r}},
\eena
where
\bena
\hs {\tilde \mu}=\hs \mu+\sqrt{\eta}\hs{\zeta}, \\
\hs \mu=(\hs{\nabla ^4\phi+2a \nabla ^2\phi+\alpha\phi)} +2\hs r\ohs g.\label {model23}\\
\eena
\end{sch}
The discrete energy is defined as follows
\bena
F^n=\int_{\Omega} [\frac{\phi^n}{2}(\nabla^4+2a\nabla^2+\alpha)\phi^{n}]\mathrm{d {\bf r}}+(r^n)^2+\frac{(\zeta^n)^2}{2}-C_0.
\eena
The scheme can be recast into
\bena
A\phi^{n+1}+(\phi^{n+1},c_1)d_1+(\phi^{n+1},c_2)d_2=b^n,
\eena
where
\bena
A=I+\Delta t \ohs M [\frac{\nabla^4}{2}+a\nabla^2+\frac{\alpha}{2}],\\
c_1=\ohs g,\\
d_1=\Delta t \ohs M \ohs g ,\\
c_2=1,\
d_2= \frac{\Delta t \ohs M}{2}\eta,\\
b^n=\phi^n-\Delta t \ohs M [\frac{\nabla^2}{2}\phi^{n}+a\nabla^2\phi^{n}+\frac{\alpha}{2}\phi^{n}+2r^n\ohs g- \\
\ohs g \int_\Omega \ohs g \phi^n \mathrm d{\bf r}]-\Delta t \ohs M\sqrt\eta\zeta^n+\Delta t \ohs M\frac{\eta}{2}\int_\Omega\phi^{n}\mathrm{d {\bf r}}.
\eena
It implies that
\bena
(\phi^{n+1},c_1)+ (\phi^{n+1},c_1) (A^{-1}d_1, c_1)+ (\phi^{n+1},c_2) (A^{-1}d_2, c_1)=(A^{-1} b^n, c_1),\\
(\phi^{n+1},c_2)+ (\phi^{n+1},c_1) (A^{-1}d_1, c_2)+ (\phi^{n+1},c_2) (A^{-1}d_2, c_2)=(A^{-1} b^n, c_2).
\eena
We solve for $(\phi^{n+1},c_1)$ and $(\phi^{n+1},c_2)$ from the above equation after we obtain
\ben
A[x,y,z]=[d_1, d_2,b^n].
\een
So, the solution is solved in the following steps,
\bena
\phi^{n+1}=z-[ (\phi^{n+1}, c_1) x+(\phi^{n+1}, c_2) y],\\
r^{n+1}=r^n+(\phi^{n+1}-\phi^n,\frac{\ohs g}{2}),\\
\zeta^{n+1}=\zeta^n+\sqrt{\eta}((\phi^{n+1},1)-(\phi^{n},1)).
\eena

\subsubsection{Numerical method for the Allen-Cahn model with a Lagrange multiplier by EQ methods}

\noindent \indent We reformulate the free energy density by introducing an intermediate variable
\bena
q=\phi^2.
\eena
Then, the free energy recast into
\bena
F
=\int_{\Omega} [\frac{\phi}{2}(\nabla^4+2a\nabla^2+\alpha)\phi+\frac{q^2}{4}]\mathrm{d {\bf r}}-L(\int_{\Omega}\phi(t)\mathrm{d {\bf r}}-\int_{\Omega}\phi(0)\mathrm{d {\bf r}}).
\eena
We rewrite (\ref{La}) as
\bena
\parl{\phi}{t}= -M \tilde \mu,\quad
\parl{q}{t}= q' \phi_t,
\eena
where \bena
\tilde \mu=\vparl{F}{\phi}=\nabla ^4\phi+2a\nabla ^2\phi+\alpha\phi+\frac{1}{2}qq'-L,\quad
L=\frac{1}{\int_{\Omega} M  \mathrm{d {\bf r}}}\int_{\Omega} M \tilde \mu\mathrm{d {\bf r}},\quad
q'=\parl{q}{\phi}.
\eena
We then discretize it using the linear modified Crank-Nicolson method in time  as follows.

\begin{sch} Given initial conditions $\phi^0,q^0$, we first compute $\phi^1, q^1$ by a first order scheme. Having computed $\phi^{n-1},q^{n-1}$, and $\phi^n,q^n$, we compute $\phi^{n+1},q^{n+1}$ as follows.
\bena
\phi^{n+1}-\phi^n= -\Delta t \overline{M}^{n+1/2} \hs{\tilde \mu},\label{AC-L-EQ}\\
q^{n+1}-q^n= \ohs{q'}(\phi^{n+1}-\phi^n).\\
\eena
where
\bena
\hs {\tilde \mu}=(\hs{\nabla ^4\phi+2a \nabla ^2\phi+\alpha\phi)} +\frac{1}{2}\hs{q} \ohs{q'}-\hs L,\\
\hs L=\frac{1}{\int_{\Omega} \ohs M  \mathrm{d {\bf r}}}\int_{\Omega} \ohs M \hs \mu\mathrm{d {\bf r}}.
\eena
\end{sch}
\noindent {\bf Remark:} $\hs L\neq \frac{L^n+L^{n+1}}{2}$.\\
Then, we have the following theorem
\begin{thm}
The mass of each phase is conserved, i.e.,
\bena
\int_{\Omega} \phi^{n+1} \mathrm{d {\bf r}}=\int_{\Omega} \phi^n \mathrm{d {\bf r}}.
\eena
\end{thm}
\noindent {\bf Proof:} Substituting the $\hs L$ into the   equation below, we have
\bena
\int_\Omega \frac{\phi^{n+1}-\phi^{n}}{\Delta t}\mathrm{d {\bf r}}\\
=\int_\Omega -\ohs M(\hs \mu-\hs L)\mathrm{d {\bf r}}=0.
\eena
This implies the masse-conservation property.

We define the discrete energy as follows
\bena
F^n=\int_{\Omega} [\frac{\phi^n}{2}(\nabla^4+2a\nabla^2+\alpha)\phi^{n}+\frac{({q}^{n})^2}{4}]\mathrm{d {\bf r}}.
\eena
The solution is solved  in the following steps
\bena
A[x,y]=[d, b^n],\\
(\phi^{n+1},c)=\frac{(y,c)}{1+(x, c)},\\
\phi^{n+1}=y-(\phi^{n+1}, c) x,\\
q^{n+1}=q^n+\ohs{q'}(\phi^{n+1}-\phi^n),\\
\eena
where
\bena
A=I+\Delta t \ohs M [\frac{\nabla^4}{2}+a\nabla^2+\frac{\alpha}{2}+\frac{1}{4}{(\ohs{q'})}^2],\\
c=\ohs{h'}\ohs{M}[\frac{\nabla^4}{2}+a\nabla^2+\frac{\alpha}{2}+\frac{1}{4}{(\ohs{q'})}^2],\\
d=-\frac{\Delta t \ohs{M}}{\int_\Omega \ohs M  \mathrm d{\bf r}},\\
b^n=\phi^n-\Delta t \ohs M (\frac{\nabla^4}{2}\phi^n+a\nabla^2\phi^n+\frac{\alpha}{2}\phi^n+\frac{1}{2}q^n\ohs {q'}-\frac{1}{4}{(\ohs{q'})}^2\phi^n-\\
\frac{\int_\Omega \ohs M(\frac{\nabla^4}{2}\phi^n+a\nabla^2\phi^n+\frac{\alpha}{2}\phi^n+\frac{1}{2}q^n\ohs {q'}-\frac{1}{4}{(\ohs{q'})}^2\phi^n)\mathrm d{\bf r} }{\int_\Omega \ohs M  \mathrm d{\bf r}})
\eena
\subsubsection{Numerical method for the Lagrangian models by SAV methods}

\noindent \indent We reformulate the free energy density by introducing an intermediate variable
\bena
r=\sqrt{\int_\Omega\frac{\phi^4}{4}\mathrm d{\bf r}+C_0}.
\eena
Then, the free energy recast into
\bena
F=\int_{\Omega} [\frac{\phi}{2}(\nabla^4+2a\nabla^2+\alpha)\phi]\mathrm{d {\bf r}}+r^2-C_0-L(\int_{\Omega}\phi(t)\mathrm{d {\bf r}}-\int_{\Omega}\phi(0)\mathrm{d {\bf r}}).
\eena
We rewrite (\ref{La}) as
\bena
\parl{\phi}{t}= -M \tilde \mu,\quad
\parl{r}{t}=\int_\Omega g \phi_t \mathrm d{\bf r},
\eena
where \bena
\tilde \mu=\vparl{F}{\phi}=\nabla ^4\phi+2a\nabla ^2\phi+\alpha\phi+2rg-L,\quad
L=\frac{1}{\int_{\Omega}  M  \mathrm{d {\bf r}}}\int_{\Omega} M \mu\mathrm{d {\bf r}},\quad
g=\parl{r}{\phi}=\frac{\phi^3}{2\sqrt{\int_\Omega\frac{\phi^4}{4}\mathrm d{\bf r}+C_0}}.
\eena
We then discretize it using the  modified Crank-Nicolson method in time  as follows.
\begin{sch} Given initial conditions $\phi^0,r^0$, we first compute $\phi^1, r^1$ by a first order scheme. Having computed $\phi^{n-1},r^{n-1}$, and $\phi^n,r^n$, we compute $\phi^{n+1},r^{n+1}$ as follows.
\bena
\phi^{n+1}-\phi^n= -\Delta t \overline{M}^{n+1/2} [(\hs{\nabla ^4\phi+2a \nabla ^2\phi+\alpha\phi)} +2\hs r \ohs g-\hs{L}],\label{EQPLAC}\\
r^{n+1}-r^n=\int_\Omega \ohs{g}(\phi^{n+1}-\phi^n)\mathrm d{\bf r}.\label{PFLAC}\\
\eena
where
\bena
\hs {\tilde \mu}=(\hs{\nabla ^4\phi+2a \nabla ^2\phi+\alpha\phi)} +2\hs r \ohs g-\hs L ,\\
\hs L=\frac{1}{\int_{\Omega} \ohs M  \mathrm{d {\bf r}}}\int_{\Omega}\ohs M \hs \mu\mathrm{d {\bf r}}.
\eena
\end{sch}
Then, we have the following theorem
\begin{thm}
The mass of each phase is conserved, i.e.,
\bena
\int_{\Omega} \phi^{n+1} \mathrm{d {\bf r}}=\int_{\Omega} \phi^n \mathrm{d {\bf r}}.
\eena
\end{thm}
\noindent {\bf Proof:} The proof is similar to that of theorem 4.1 and is thus omitted.

We define the discrete energy as follows
\bena
F^n=\int_{\Omega} [\frac{\phi^n}{2}(\nabla^4+2a\nabla^2+\alpha)\phi^{n}]\mathrm{d {\bf r}}+({r}^{n})^2-C_0.
\eena
This scheme can be recast into
\ben
A\phi^{n+1}+(\phi^{n+1}, c_1)d_1+(\phi^{n+1},c_2)d_2+(c_3,(\phi^{n+1},c_1))d_2=b^n,
\een
So we have
\bena
(\phi^{n+1},c_1)+(\phi^{n+1},c_1)(A^{-1} d_1, c_1)+(\phi^{n+1},c_2)(A^{-1} d_2, c_1)+(c_3,(\phi^{n+1},c_1))(A^{-1}d_2,c_1)=(A^{-1} b^n,c_1),\\
(\phi^{n+1},c_2)+(\phi^{n+1},c_1)(A^{-1} d_1, c_2)+(\phi^{n+1},c_2)(A^{-1} d_2, c_2)+(c_3,(\phi^{n+1},c_1))(A^{-1}d_2,c_2)=(A^{-1} b^n,c_2).\\
\eena
We solve for $(\phi^{n+1},c_1)$ and $(\phi^{n+1},c_2)$ from the above equations after we obtain
\bena
A[x,y,z]=[d_1,d_2, b^n],\\
\eena
where
\bena
A=I+\Delta t \ohs M(\frac{\nabla^4}{2} +a \nabla^2+\frac{\alpha}{2}),\\
c_1=\ohs g,
d_1=\Delta t \ohs M \ohs g,\\c_2=\ohs {h'}\ohs M(\frac{\nabla^4}{2} +a \nabla^2+\frac{\alpha}{2}),\\
d_2=-\Delta t \ohs M\frac{1}{\int_\Omega \ohs{M}\mathrm d{\bf r}},\\
c_3=\ohs M ,\\
b^n=\phi^n-\Delta t\ohs M(\frac{\nabla^4}{2}\phi^n +a \nabla^2\phi^n+\frac{\alpha}{2}\phi^n+2r^n\ohs{g}-\ohs{g}\int_\Omega\ohs{g}\phi^n\mathrm{d {\bf r}})\\
-(\ohs M,\frac{\nabla^4}{2}\phi^n +a \nabla^2\phi^n+\frac{\alpha}{2}\phi^n+2r^n\ohs{g}-\ohs{g}\int_\Omega\ohs{g}\phi^n\mathrm{d {\bf r}})d_2.\\
\eena
Finally, the solution is obtained as follows
\bena
A[x,y,z]=[d_1,d_2, b^n],\\
\phi^{n+1}=z-(\phi^{n+1}, c_1)x-(\phi^{n+1},c_2)y-(c_3,(\phi^{n+1},c_1))y,\\
r^{n+1}=r^n+(\ohs g,(\phi^{n+1}-\phi^n)).
\eena

\subsection{Spatial discretization}

\noindent \indent We use  the finite difference method to discretize the equations with the Neumann boundary condition in space for all models.
The linear spatially dependent PDE systems resulting from all semi-discrete schemes are spatially discretized by compact second order finite difference methods at the cell center. We divide the 2D domain $\Omega=[0,L_x] \times [0, L_y]$ into rectangular meshes with mesh sizes $h_x=L_x/N_x$ and $h_y=L_y/N_y$, where $L_x$, $L_y$ are two positive real numbers and $N_x$, $N_y$ are the number of meshes in each direction. After this, the sets of the cell center points $C_x$ and $C_y$ for the uniform partition are defined as follows
\bena
C_x=\left\{x_{i}|i=0, 1, \cdots, N_x \right\},\quad
C_y=\left\{y_{j}|j=0, 1, \cdots, N_y \right\},
\eena
where $x_i=(i-\frac{1}{2})h_x$ and $y_j=(j-\frac{1}{2})h_y$.

We define the east-west-edge-to-center and center-to-east-west-edge difference operators $d_x$ and $D_x$ as follows
\bena
d_x\phi_{ij}=\frac{\phi_{i+\frac{1}{2},j}-\phi_{i-\frac{1}{2},j}}{h_x}, \qquad D_x\phi_{i+\frac{1}{2,}j}=\frac{\phi_{i+1,j}-\phi_{i-1,j}}{h_x}.
\eena
Similarly, we define the north-south-edge-to-center and center-to-north-south-edge difference operators $d_y$ and $D_y$ as follows
\bena
d_y\phi_{ij}=\frac{\phi_{i,j+\frac{1}{2}}-\phi_{i,j-\frac{1}{2}}}{h_y}, \qquad D_y\phi_{i,j+\frac{1}{2}}=\frac{\phi_{i,j+1}-\phi_{i,j-1}}{h_y}.
\eena
The fully discrete Laplacian and fourth order gradient operator are given by
\bena
\nabla_h^2=\Delta_h=d_x(D_x\phi)+d_y(D_y\phi),
\nabla_h^4=\Delta_h^2=d_x(D_x\Delta_h)+d_y(D_y\Delta_h).
\eena
The discrete inner product   is defined as follows
\bena
<f,g>=h_xh_y\sum_{i,j}f_{i,j}g_{i,j},\label{DiscreteInnerproduct}
\eena
where $f_{ij}$ and $g_{ij}$ are given at the cell center.
In particular,
\bena
<f,1>=h_xh_y\sum_{i,j}f_{i,j},\quad
\lVert f \rVert_d=\sqrt{<f,f>}.
\eena
With the notations, we replace the continuous differential operators in the semi-discrete schemes by the discrete difference operators to arrive at the corresponding fully discrete schemes. To save space, we will not enumerate them here.

\subsection{Energy dissipation property and solvability of the linear systems resulting from the schemes}

\noindent \indent We summarize the energy dissipation law and unique solvability for all linear systems resulting from the semidiscrete schemes presented in this section into two theorems.
The proofs of the energy dissipation property and solvability for the fully discrete schemes are similar, we  only prove the theorems  for the Allen-Cahn model  with a Lagrangian multiplier discreized using EQ methods below and omit others for simplicity. Based on scheme 4.7 and the spatial discretization, the fully discrete scheme corresponding to scheme 4.7 is summarized below.

\begin{sch} Given initial conditions $\phi^0,q^0$, we first compute $\phi^1, q^1$ by a first order scheme. Having computed $\phi^{n-1},q^{n-1}$, and $\phi^n,q^n$, we compute $\phi^{n+1},q^{n+1}$ as follows.
\bena
\phi^{n+1}-\phi^n  -\Delta t \overline{M}^{n+1/2} [(\hs{\nabla_h ^4\phi+2a \nabla_h ^2\phi+\alpha\phi)} +\frac{1}{2}\hs{q} \ohs{q'}-\hs{L}],\label{AC-L-EQ-F}\\
q^{n+1}-q^n= \ohs{q'}(\phi^{n+1}-\phi^n).\\
\eena
where
\bena
\hs {\tilde \mu}=(\hs{\nabla_h ^4\phi+2a \nabla_h ^2\phi+\alpha\phi)} +\frac{1}{2}\hs{q} \ohs{q'}-\hs L,\\
\hs L=\frac{1}{\int_{\Omega} \ohs M  \mathrm{d {\bf r}}}\int_{\Omega} \ohs M \hs \mu\mathrm{d {\bf r}}.
\eena
\end{sch}
We define the discrete energy as follows
\bena
F^n=<\frac{\phi^n}{2}(\nabla_h^4+2a\nabla_h^2+\alpha)\phi^{n}+\frac{({q}^{n})^2}{4},1>.
\eena
\begin{thm}
The fully discrete  scheme obeys the following energy dissipation law
\bena
F^{n+1}-F^n=-\Delta t<\hs \mu-\hs L,\ohs{M}(\hs \mu-\hs L )>.
\eena
\end{thm}
\noindent {\bf Proof:}
Taking inner product of (\ref{AC-L-EQ-F}) with $-u^{n+1/2}$, we obtain
\bena
-<\frac{\phi^{n+1}-\phi^n}{\Delta t},\mu^{n+1/2}>\\
=-<-\ohs{M}[\hs{\mu}-\hs{L}],\hs{\mu}>\\
=\lVert\sqrt{\ohs{M}}(\hs{\mu}-\hs{L})\rVert_d^2+<\ohs{M}[\hs{\mu}-\hs{L}],\hs{L}>.
\eena
Taking inner product of $\mu^{n+1/2}$ with $ \frac{\phi^{n+1}-\phi^n}{\Delta t} $, we have
\bena
<(\hs{\nabla_h ^4\phi+2a \nabla_h ^2\phi+\alpha\phi)}+\frac{1}{2}\hs{q} \ohs{q'},\frac{\phi^{n+1}-\phi^n}{\Delta t}>\\
=\frac{\lVert\nabla_h ^2\phi^{n+1}\rVert_d^2-\lVert\nabla_h ^2\phi^n\rVert_d^2}{2\Delta t}-a\frac{\lVert\nabla_h \phi^{n+1}\rVert_d^2-\lVert\nabla_h \phi^n\rVert_d^2}{\Delta t}+\frac{\alpha}{2 \Delta t}(\lVert\phi^{n+1}\rVert_d^2-\lVert\phi^n\rVert_d^2)+<\frac{1}{2}\hs{q} \ohs{q'},\frac{\phi^{n+1}-\phi^n}{\Delta t}>.
\eena
Taking  inner product of $q^{n+1}-q^{n}$ with $ \frac{q^{n+1}+q^n}{\Delta t} $, we obtain
\bena
\frac{1}{\Delta t}(\lVert q^{n+1}\rVert_d^2-\lVert q^{n}\rVert_d^2)=\frac{1}{\Delta t}<\ohs{q'}(\phi^{n+1}-\phi^n),q^{n+1}+q^{n}>.
\eena
Combining the above equations, we obtain
\bena
\frac{\lVert\nabla_h ^2\phi^{n+1}\rVert_d^2-\lVert\nabla_h ^2\phi^n\rVert_d^2}{2\Delta t}-a\frac{\lVert\nabla_h \phi^{n+1}\rVert_d^2-\lVert\nabla_h \phi^n\rVert_d^2}{\Delta t}+\frac{\alpha}{2 \Delta t}(\lVert\phi^{n+1}\rVert_d^2-\lVert\phi^n\rVert_d^2)+\frac{1}{4\Delta t}(\lVert q^{n+1}\rVert_d^2-\lVert q^{n}\rVert_d^2)\\
=-\lVert\sqrt{\ohs{M}}(\hs{\mu}-\hs{L})\rVert_d^2-(\ohs{M}[\hs{\mu}-\hs{L}],\hs{L}).
\eena
Substituting the expression of \hs L into the equation, we have
\bena
-\lVert\sqrt{\ohs{M}}(\hs{\mu}-\hs{L})\rVert_d^2-(\ohs{M}[\hs{\mu}-\hs{L}],\hs{L})\\
=-\lVert\sqrt{\ohs{M}}(\hs{\mu}-\hs{L})\rVert_d^2.
\eena
\noindent {\bf Remark:} (i). This  proof applies to  the semi-discrete schemes as well.
(ii). When the linear schemes of the nonlocal Allen-Cahn model involve  integrals discretized by a composite Trapezoidal rule, efficient numerical methods can be devised to solve the resulting linear systems. Such methods are derived from the Sherman-Morisson formula (See Appendix).

\begin{thm}
The linear system resulting from the above fully discrete scheme admits a unique solution.
\end{thm}
\noindent {\bf Proof:}
Note that the solution in scheme \ref{AC-L-EQ-F} is solved  via the following steps
\bena
A\phi^{n+1}+(\phi^{n+1},c)d=b^n,\quad
q^{n+1}-q^n= \ohs{q'}(\phi^{n+1}-\phi^n),
\eena
where
\bena
A=I+\Delta t \ohs M [\frac{\nabla_h^4}{2}+a\nabla_h^2+\frac{\alpha}{2}+\frac{1}{4}{(\ohs{q'})}^2],\\
c=\ohs{M}[\frac{\nabla_h^4}{2}+a\nabla_h^2+\frac{\alpha}{2}+\frac{1}{4}{(\ohs{q'})}^2],\\
d=-\frac{\Delta t \ohs{M}}{<\sqrt{\ohs M},\sqrt{\ohs M>}},\\
b^n=\phi^n-\Delta t \ohs M (\frac{\nabla_h^4}{2}\phi^n+a\nabla_h^2\phi^n+\frac{\alpha}{2}\phi^n+\frac{q^n\ohs {q'}}{2}-\frac{1}{4}{(\ohs{q'})}^2\phi^n-\\
\frac{<\ohs M,\frac{\nabla_h^4}{2}\phi^n+a\nabla_h^2\phi^n+\frac{\alpha}{2}\phi^n+\frac{q^n\ohs {q'}}{2}-\frac{1}{4}{(\ohs{q'})}^2\phi^n> }{<\sqrt{\ohs M},\sqrt{\ohs M}>}).
\eena
From the Sherman-Morrison formula, we notice that the solution uniqueness of $A\phi^{n+1}+(\phi^{n+1},c)d=b^n$ depends on the uniqueness of the  corresponding linear system $A\phi^{n+1}=b^n$. Now we only need to prove the uniqueness of the solution for
\bena
A\phi^{n+1}={\bf 0}.
\eena
 If $\Delta t$ is small enough, we have
\bena
<A\phi,\phi>=<\phi,A\phi>=<\phi,(I+\Delta t \overline M [\frac{\nabla_h^4}{2}+a\nabla_h^2+\frac{\alpha}{2}+\frac{1}{4}{(\overline{q'})}^2])\phi>\\=<\phi,\phi>+<\phi,\Delta t \overline M \frac{\nabla_h^4}{2}\phi>
+<\phi,\Delta t \overline M a\nabla_h^2\phi>\\
+<\phi,\Delta t \overline M \frac{\alpha}{2}\phi>
+<\phi,\Delta t \overline M \frac{1}{4}{(\overline{q'})}^2\phi>\\
\geq  0.
\eena
So, $A\phi=0$ has only zero solution.

\noindent {\bf Remark:} One can prove the uniqueness of the solution for any time step size  if $\hs \mu=\nabla_h ^4 \hs \phi+2a \nabla_h ^2 \ohs\phi+\alpha \hs \phi +\frac{1}{2}\hs{q} \ohs{q'}$ in equation \eqref{AC-L-EQ-F}.

\section{Numerical Results and Discussions}

\noindent \indent In this section, we conduct  mesh refinement tests  to validate the accuracy of the proposed schemes and then present some numerical examples to assess the schemes for the nonlocal Allen-Cahn models against those for the Cahn-Hilliard model. For convenience, we refer the numerical schemes designed by EQ methods for the Allen-Cahn model, the Cahn-Hilliard model, the  Allen-Cahn model with a penalizing potential and the  Allen-Cahn model with a Lagrangian multiplier as AC-EQ, CH-EQ, AC-P-EQ and AC-L-EQ, respectively. Similarly, we name the numerical schemes obtained using SAV approaches for the models as AC-SAV, CH-SAV, AC-P-SAV, AC-L-SAV, respectively. In the following, we set the constant in the free energy at $C_0=1\times 10^ 4$ in all computations.

\subsection{Accuracy test}

\noindent \indent We confirm the convergence rates of the proposed schemes for the PFC models through mesh refinement tests.
The computational domain is set as $\Omega=[0,1]^2$ . The model parameter values are chosen as $a=1, \varepsilon=0.1, M=1\times 10^{-3}$.
We solve the equations with the initial condition given by
\bena
\phi(0,x,y)=\frac{1}{2}+\frac{1}{2}\cos(\pi x)\cos(\pi y).
\eena
We choose the space step size $h_x=h_y=\frac{1}{256}$. By taking a linear refinement path $\Delta t=\frac{0.05}{2^k}$, $k=0, 1 , \cdots, 5$,  we calculate the $L^2$ errors of the phase variable with adjacent k at $t=1$. The tables show the schemes are second order accurate in time numerically.

We also compare the computational efficiency of all schemes designed by EQ and SAV methods with $M=1\times 10^{-3},1$ and $1\times 10^{2}$ in table 3. The AC-P-EQ/SAV schemes perform the best among the schemes for nonlocal Allen-Cahn models. Besides this, AC-P-EQ/SAV schemes also perform better than CH-EQ/SAV schemes in most test cases. In fact, the accuracy of the Cahn-Hillard model relies on the mobility coefficient $M$ more sensitively than the Allen-Cahn models do. Hence, the accuracy of the schemes for the nonlocal Allen-Cahn models is better than that for the Cahn-Hillard model if $M$ is large. We will  discuss it in more details next (see Figure \ref{Fig6}).
\vskip 10 pt
{
\centerline{{\bf Table 1} Mesh refinement tests for the proposed schemes using EQ methods.}
\vskip 10 pt
\centerline{
\begin {tabular} [!htbp] {|c|c|c|c|c|c|c|c|c|c|} 
\hline
 \multicolumn {2}{|c|} {Scheme}& \multicolumn {2}{|c|} {AC-EQ} & \multicolumn {2}{|c|} {CH-EQ} & \multicolumn {2}{|c|} {AC-P-EQ}& \multicolumn {2}{|c|} {AC-L-EQ}\\
\cline{0-9}
\hline
Coarse $\Delta t$&Fine $\Delta t$ &$L^2$ error & order& $L^2$ error & order & $L^2$ error & order& $L^2$ error & order\\
\cline{0-9}
5.00E-02&2.5E-2&1.21E-06&-&1.23E-05 &-&1.21E-06 &- &1.21E-06&- \\
2.5E-2&1.25E-2&3.02E-07&2.00&3.12E-06 &1.98&3.02E-07 &2.00 &3.02E-07&2.00  \\
1.25E-2&6.25E-3&7.56E-08&2.00&7.84E-07 &1.99&7.56E-08 &2.00 &7.56E-08&2.00  \\
6.25E-3&3.125E-3&1.89E-08&2.00&1.96E-07 &2.00&1.89E-08 &2.00 &1.89E-08&2.00  \\
3.125E-3&1.5625E-3&4.72E-09&2.00&4.90E-08 &2.00&4.72E-09 &2.00 &4.72E-09&2.00 \\
\hline
\end {tabular}
}}
\vskip 10 pt
{
\centerline{{\bf Table 2} Mesh refinement tests for the proposed schemes using SAV methods.}
\vskip 10 pt
\centerline{
\begin {tabular} [!htbp] {|c|c|c|c|c|c|c|c|c|c|} 
\hline
 \multicolumn {2}{|c|} {Scheme}& \multicolumn {2}{|c|} {AC-SAV} & \multicolumn {2}{|c|} {CH-SAV} & \multicolumn {2}{|c|} {AC-P-SAV}& \multicolumn {2}{|c|} {AC-L-SAV}\\
\cline{0-9}
\hline
Coarse $\Delta t$&Fine $\Delta t$ &$L^2$ error & order& $L^2$ error & order & $L^2$ error & order& $L^2$ error & order\\
\cline{0-9}
5.00E-02&2.5E-2&1.24E-06&-&1.22E-05 &-&1.24E-06&- &1.24E-06&- \\
2.5E-2&1.25E-2&3.09E-07&2.00&3.10E-06 &1.98&3.09E-07 &2.00 &3.09E-07&2.00  \\
1.25E-2&6.25E-3&7.72E-08&2.00&7.78E-07 &1.99&7.72E-08&2.00 &7.72E-08&2.00  \\
6.25E-3&3.125E-3&1.93E-08&2.00&1.95E-07 &2.00&1.93E-08&2.00 &1.93E-08&2.00  \\
3.125E-3&1.5625E-3&4.83E-09&2.00&4.87E-08 &2.00&4.82E-09&2.00 &4.82E-09&2.00 \\
\hline
\end {tabular}
}}
\vskip 10 pt
\centerline{{\bf Table 3} Efficiency of the schemes at $1000$ time steps with respect to  $M=0.001,1$ and $100$ (From top to below).}
\vskip 10 pt
\centerline{
\begin {tabular} [!htbp] {|c|c|c|c|c|c|c|c|c|} 
\hline
Scheme &  AC-EQ &AC-SAV& CH-EQ& CH-SAV &  AC-P-EQ &AC-P-SAV & AC-L-EQ&AC-L-SAV\\
\cline{0-8}
\hline
\cline{0-8}
\hline
Time (s)& 45.3&77.9&66.0&107.8&68.3&90.0&110.0&103.2\\
\cline{0-8}
\hline
Time (s)& 44.4&76&54.6&104&38.8&84.2&53.9&100.7\\
\cline{0-8}
Time (s)& 53.4&75.0&66.8&111.9&61.4&83.2&103.3&99.5\\
\hline				
\end {tabular}
}
\subsection{Assessment of the numerical schemes }
\noindent \indent To further assess the numerical schemes, we numerically solve the model equations using the schemes with respect to  two benchmark problems.
 Firstly we simulate the phase transition of crystal growth in 2D. We use time step $\Delta t=1 \times 10^{-3}$ and  $256\times 256$ space mesh in the 2D simulation.  A solid crystallite with Hexagonal ordering in 2D is initially placed in the centre of the domain, which is assigned an average density $\overline{\phi}$. The initial condition is given by
\bena
\phi_0(\bf r)=\overline{\phi}+w(x)(A\phi_s(\bf r)),
\eena
where
\bena
w(\bf r)=\begin{cases} (1-(\frac{|\bf r-\bf r_0| }{\bf d_0})^2)^2& \text{if }  \frac{|\bf r-\bf r_0|}{\bf d_0}\le 1, \\
0& \text{otherwise}. \end{cases}
\eena
\bena
\phi_s (\bf r)=\cos(\frac{q}{\sqrt{3}}y)\cos(qx)-\frac{1}{2}\cos(\frac{2q}{\sqrt 3}y),
\eena
$\bf r_0$ is the center coordinate of  the domain, and $\bf d_0$ is $\frac{1}{6}$ of the domain length in the x-direction. The domain is given by  $ \Omega =[0,\frac{2\pi}{q}a] \times [0,\frac {\sqrt {3} \pi}{q}b]$, $a=10$ and $b=12$. The other values are  $\varepsilon=0.325, \overline{\phi}=\frac{\sqrt {\varepsilon}}{2}, A=\frac{4}{5}(\overline{\phi}+\frac{\sqrt{15\varepsilon-36\overline {\phi}^2}}{3})$ and $q=\frac{\sqrt{3}}{2}.$

\begin{figure*}
\centering
\subfigure[]{
\begin{minipage}[b]{0.49\linewidth}
\includegraphics[width=1\linewidth]{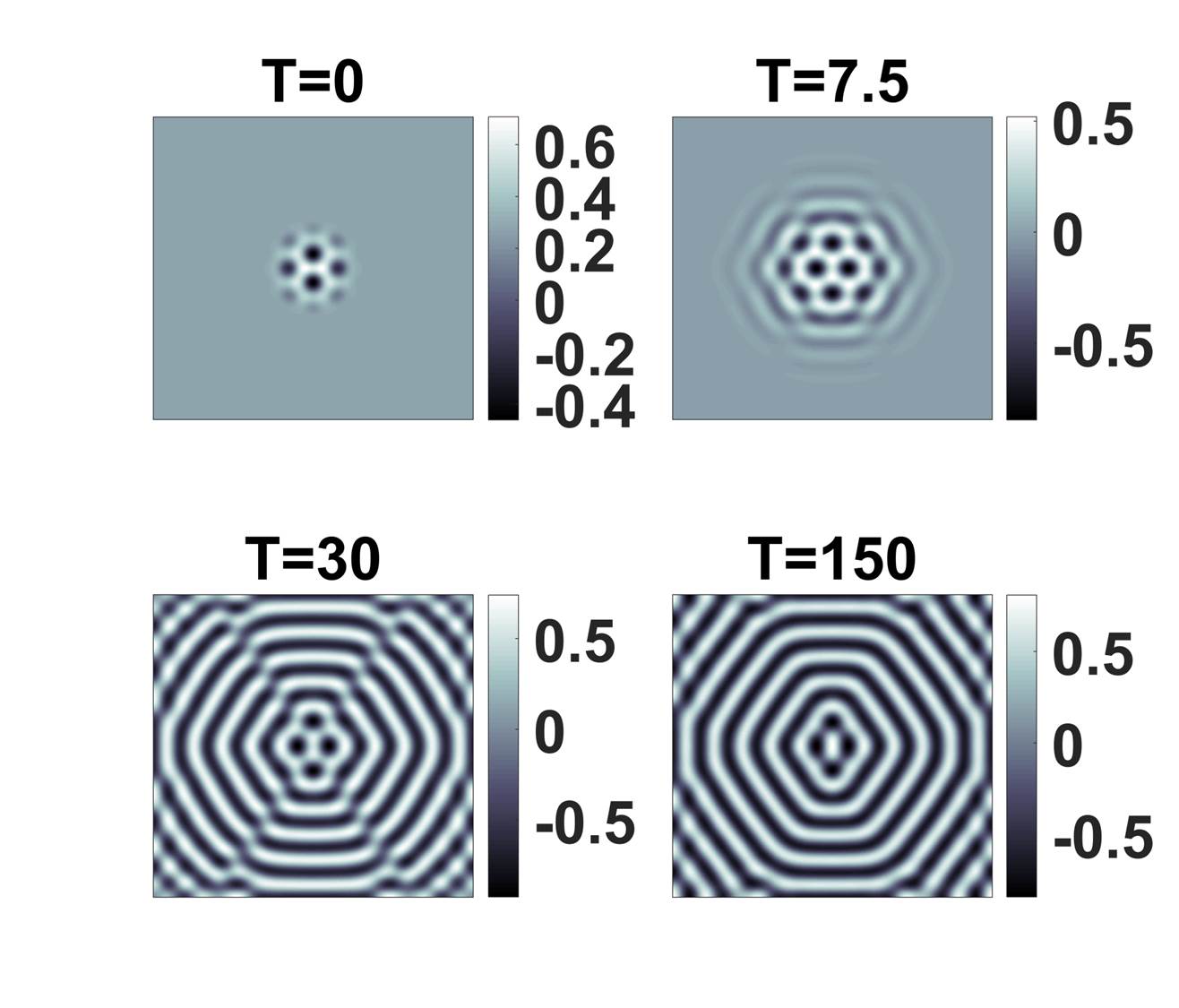}
\end{minipage}}
\subfigure[]{
\begin{minipage}[b]{0.49\linewidth}
\includegraphics[width=1\linewidth]{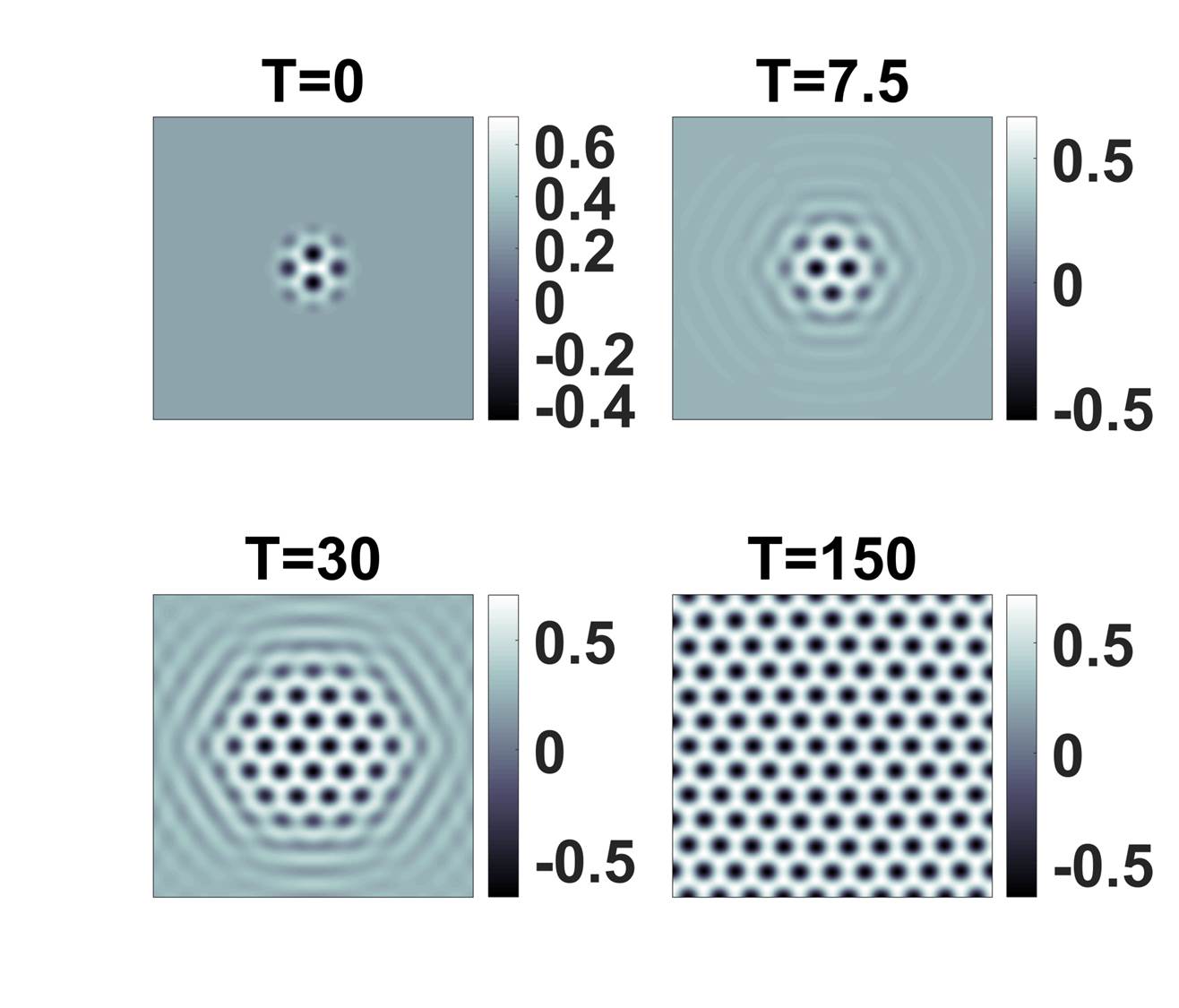}
\end{minipage}}
\subfigure[]{
\begin{minipage}[b]{0.49\linewidth}
\includegraphics[width=1\linewidth]{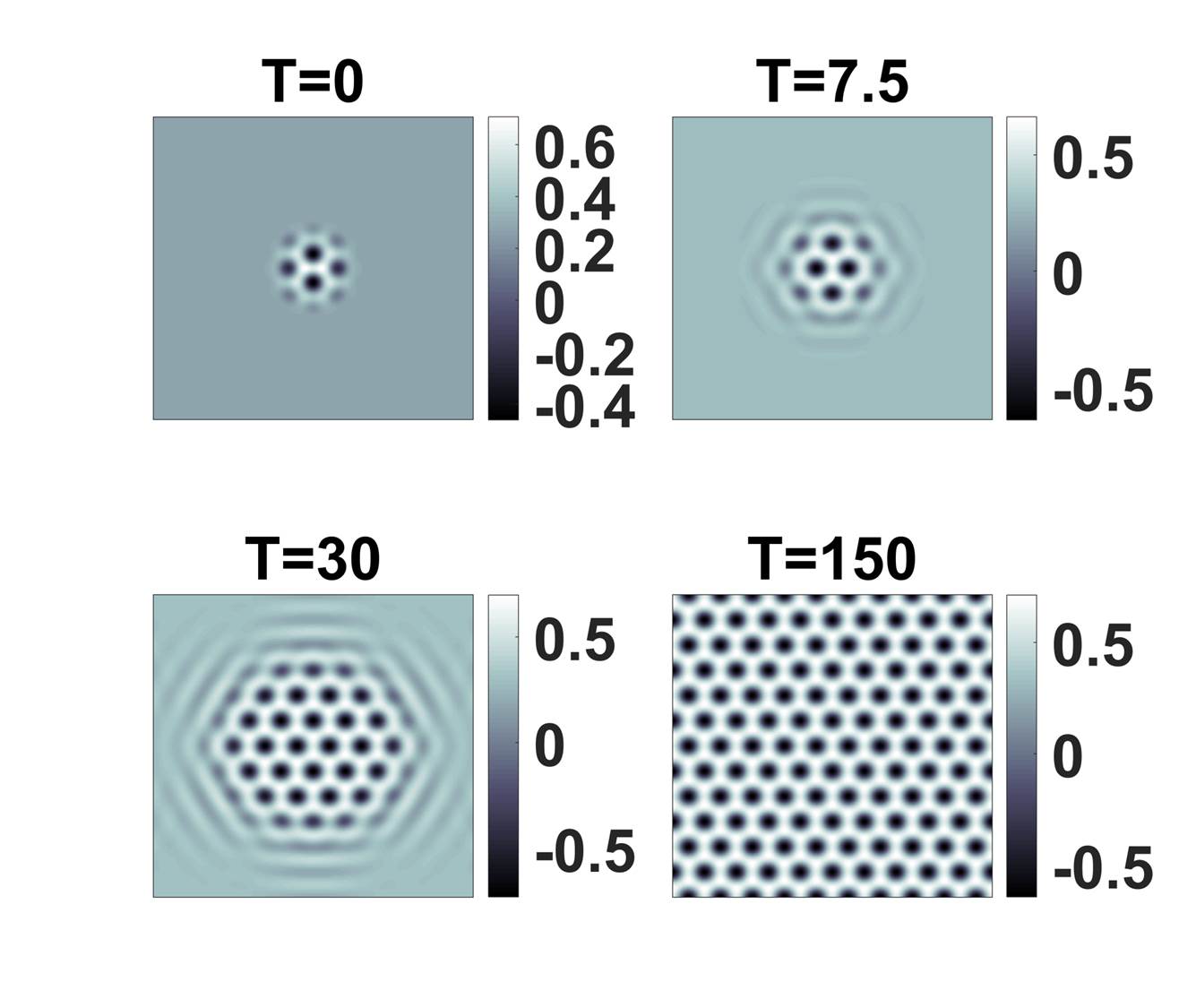}
\end{minipage}}
\subfigure[]{
\begin{minipage}[b]{0.49\linewidth}
\includegraphics[width=1\linewidth]{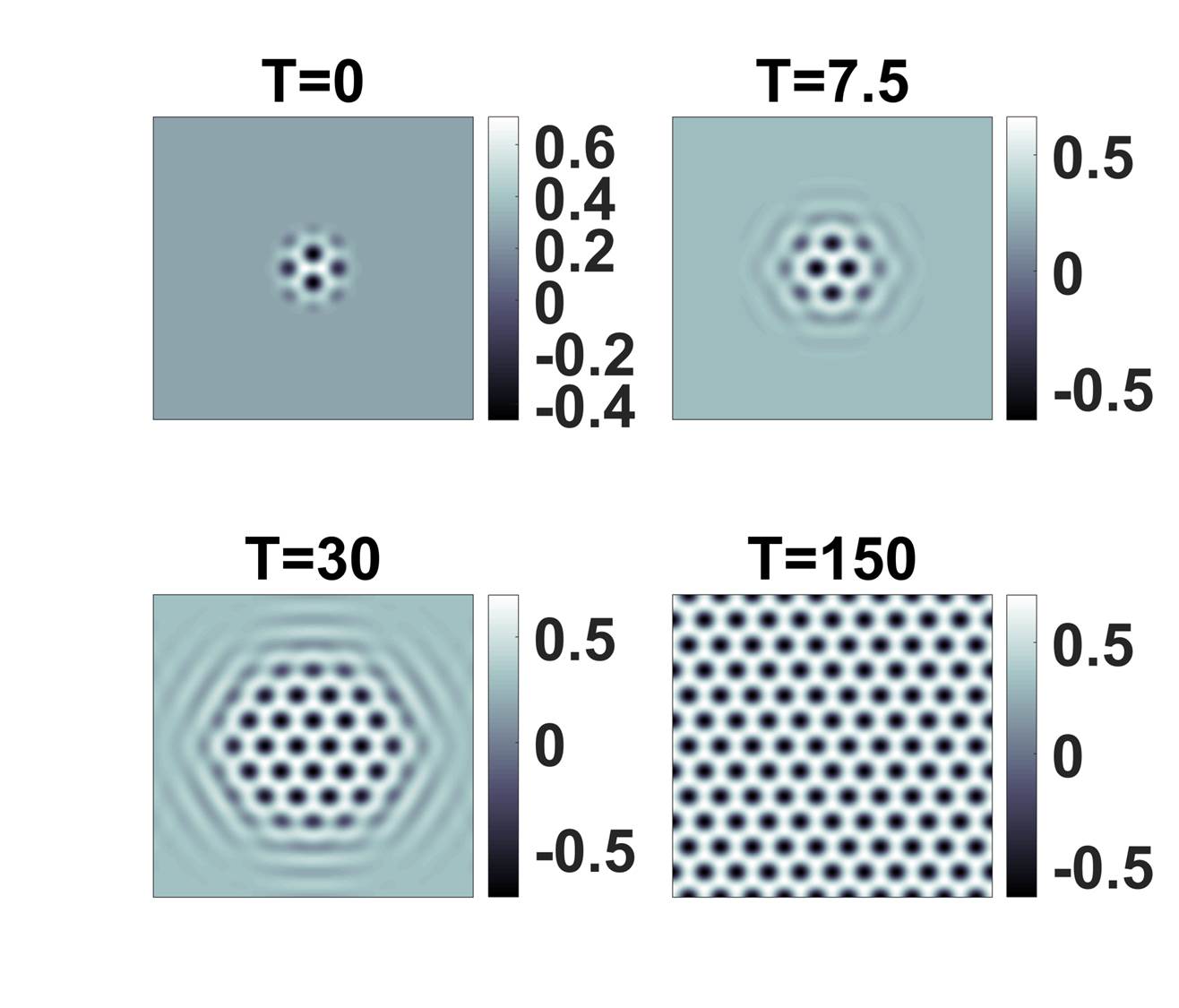}
\end{minipage}}
\caption{Crystal growth simulated using the Allen-Cahn, Cahn-Hillard and Allen-Cahn models with nonlocal constraints at $M=1$.  (a)-(d) are computed using AC-EQ, CH-EQ, AC-P-EQ and AC-L-EQ schemes respectively. Snapshots of the  atomistic density field $\phi$ are depicted at $\mathrm T=0, 7.5, 30, 150$, respectively. Parameter $\eta$ is set as $1\times 10^3$ in the AC-P-EQ/SAV schemes. We use time step $\Delta t=1 \times 10^{-3}$ and  $256\times 256$ space meshes in the 2D simulation. The Allen-Cahn model gives an erroneous result while the other models give comparable results.}\label{Fig1}
\end{figure*}

\begin{figure*}
\centering
\subfigure[]{
\begin{minipage}[b]{0.49\linewidth}
\includegraphics[width=1\linewidth]{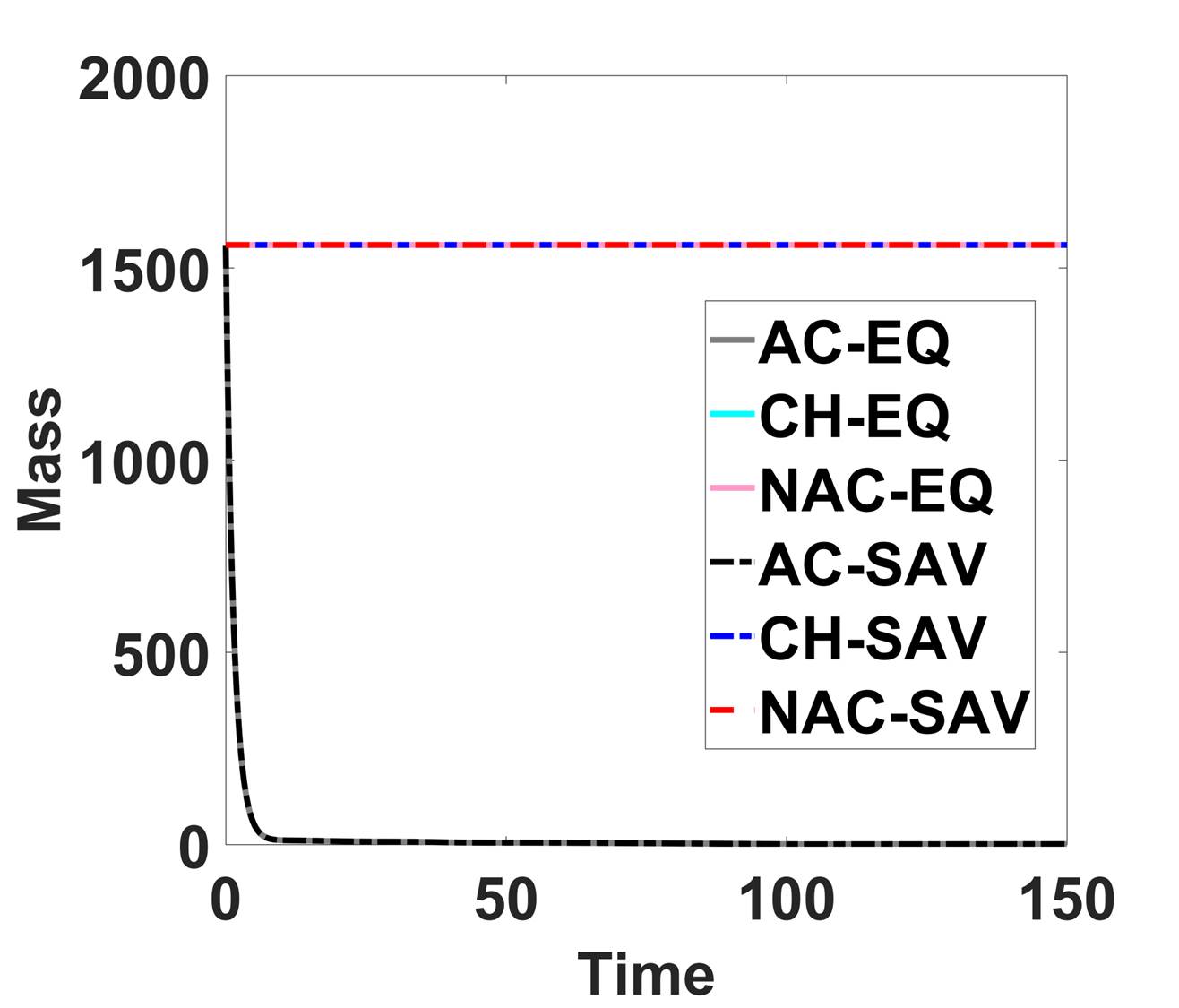}
\end{minipage}}
\subfigure[]{
\begin{minipage}[b]{0.49\linewidth}
\includegraphics[width=1\linewidth]{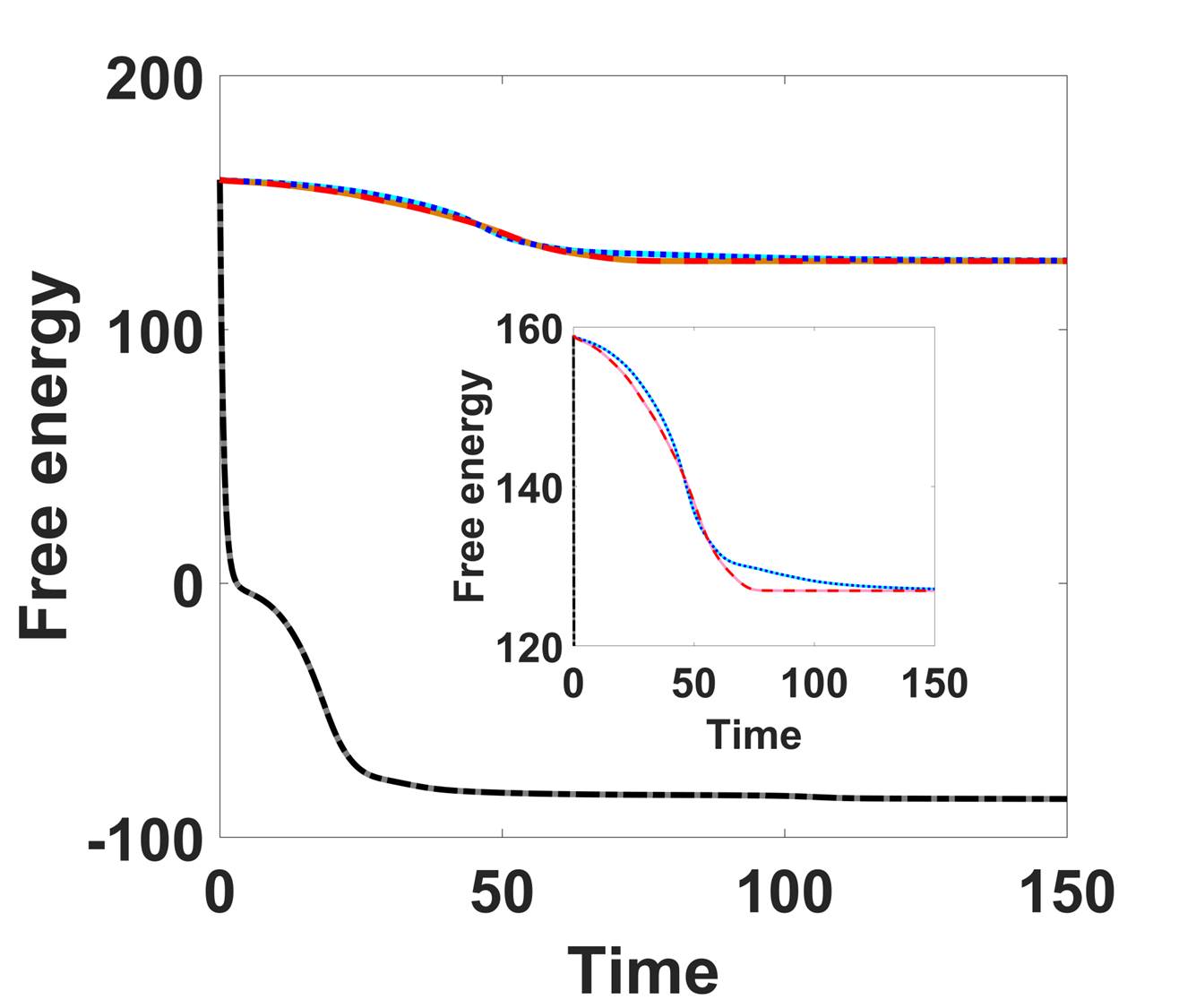}
\end{minipage}}
\caption{Time evolution of the mass and free energy from the simulations for Allen-Cahn (AC), Cahn-Hillard (CH) and nonlocal Allen-Cahn (NAC) models are shown in (a) and (b), respectively. Since the results computed by the Allen-Cahn model with nonlocal constraints are nearly identical, we only show time evolution of the mass and free energy computed using the AC-L-EQ scheme. We compare time evolution of the mass and free energy computed by the EQ and SAV schemes in (a) and (b). There is no difference between the results computed by the two methods. The mass computed using the Allen model vanishes before $\mathrm T=10$, whereas is conserved in the Cahn-Hillard and the Allen-Cahn model with nonlocal constraints. The free energy computed by all models are dissipative. The nonlocal Allen-Cahn models predict  comparable   free energy to the Cahn-Hillard model but reach the steady state faster than the Cahn-Hillard model does.}\label{Fig2}
\end{figure*}

 Figure \ref{Fig1} shows time evolution  of the crystal growth process computed by AC-EQ, CH-EQ, AC-P-EQ and AC-L-EQ schemes, respectively. In Figure \ref{Fig1}-(a), the crystal growth simulated by the Allen-Cahn model can't preserve the Hexagonal ordering, different from the results simulated by the Cahn-Hillard model in Figure \ref{Fig1}-(b) and the nonlocal Allen-Cahn models in  Figure \ref{Fig1}-(c) and (d). The time evolution of mass and free energy are shown in Figure \ref{Fig2}-(a) and (b) respectively. The results computed by the EQ and SAV schemes for the same  model (Allen-Cahn model, Cahn-Hillard model and nonlocal Allen-Cahn model) are identical. We don't see any differences between the results of the Allen-Cahn model with a penalizing potential and the Allen-Cahn model with a Lagrangian multiplier either. The mass   decays in the Allen-Cahn model   to nearly zero in finite time. In contrast, the mass in the Cahn-Hillard model and the nonlocal Allen-Cahn models is conserved in the simulations. Meantime, the free energies of the Cahn-Hillard model and the nonlocal Allen-Cahn models are larger than that of the Allen-Cahn model. Figure \ref{Fig2}-(b) shows that the free energy computed by the nonlocal Allen-Cahn model   reaches the steady state faster than that of the Cahn-Hillard model.

 Secondly we simulate another case of polycrystalline growth involving the grain boundary effect, where the two initial crystallites with a hexagonal configuration oriented in different direction (or misorientation) are put in the domain. Grain boundaries appear when the two crystals meet during the growth, which yields some orientation mismatch. The initial condition is given by
\bena
\phi_0(\bf r)=\overline{\phi}+w(x)(A\phi_s(\bf r)),
\eena
where
\bena
w(\bf r)=\begin{cases} (1-(\frac{|\bf r_1| }{\bf d_0})^2)^2& \text{if }  \frac{|\bf r_1|}{\bf d_0}\le 1, \\
(1-(\frac{|\bf r_2| }{\bf d_0})^2)^2& \text{if }  \frac{|\bf r_2|}{\bf d_0}\le 1, \\
0& \text{otherwise}, \end{cases}
\eena
${\bf r_1}=\sqrt{(x-\frac{1}{2} x_0)^2+(y-\frac{1}{2}y_0)^2}$, ${\bf r_2}=\sqrt{(x-\frac{3}{2} x_0)^2+(y-\frac{3}{2}y_0)^2}$ and $(x_0,y_0)$ is the center of the domain. The other parameters and the domain are the same as in the first example.
By doing an affine transformation of the Cartesian coordinates $(x,y)$ to produce a rotation $\theta$ in the domain, the modified coordinates $(x_{\theta},y_{\theta})$ can be used to generate the crystallites in different directions,
\bena
x_{\theta}=cos(\theta)x-sin(\theta)y,\\
y_{\theta}=sin(\theta)x+cos (\theta)y. \label{Orientation}
\eena
We put two crystallites in the domain, the first one is defined as equation (\ref{Orientation}) with $\theta=0$, the other is with $\theta=\frac{\pi}{8}.$

\begin{figure*}
\centering
\subfigure[]{
\begin{minipage}[b]{0.49\linewidth}
\includegraphics[width=1\linewidth]{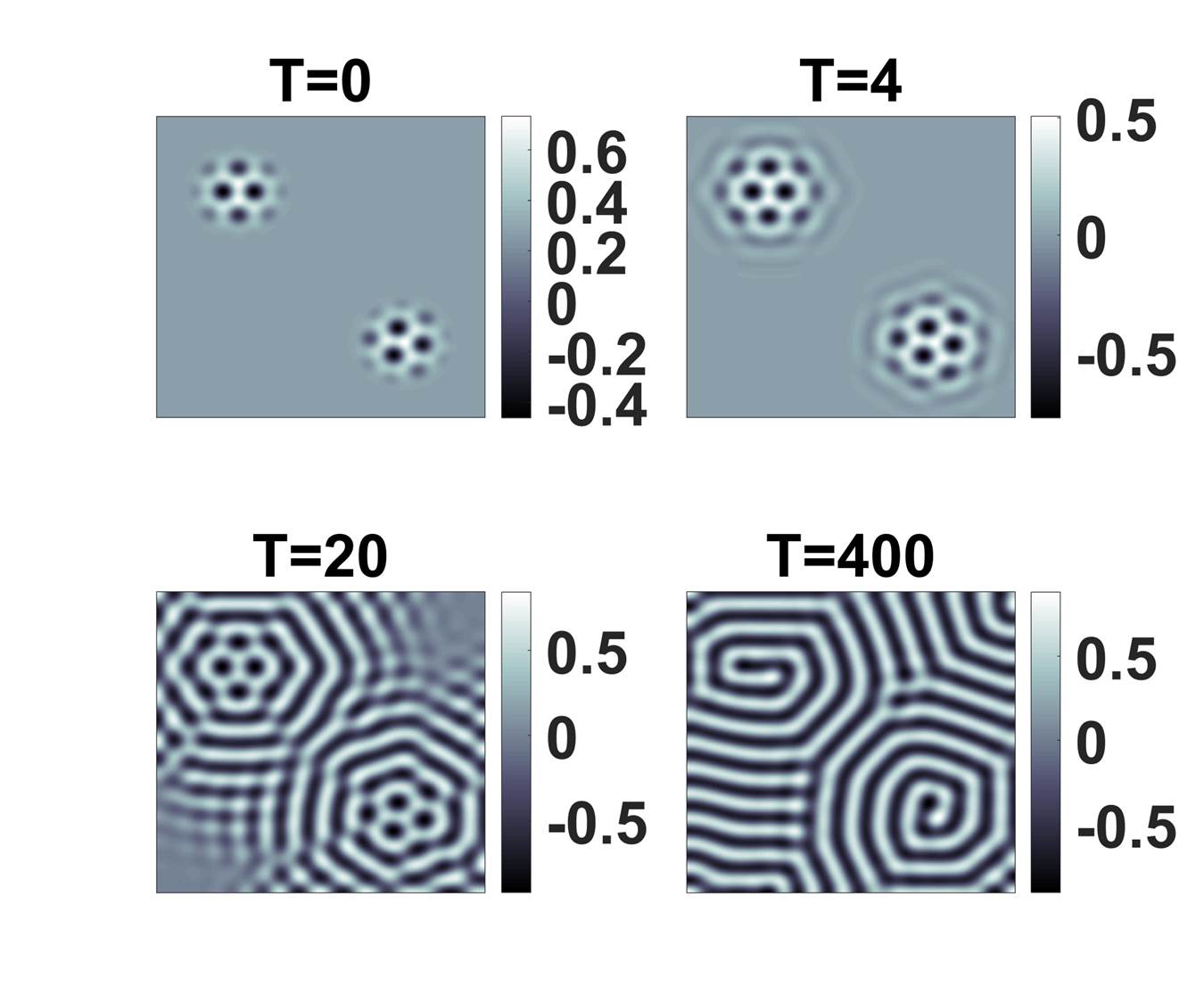}
\end{minipage}}
\subfigure[]{
\begin{minipage}[b]{0.49\linewidth}
\includegraphics[width=1\linewidth]{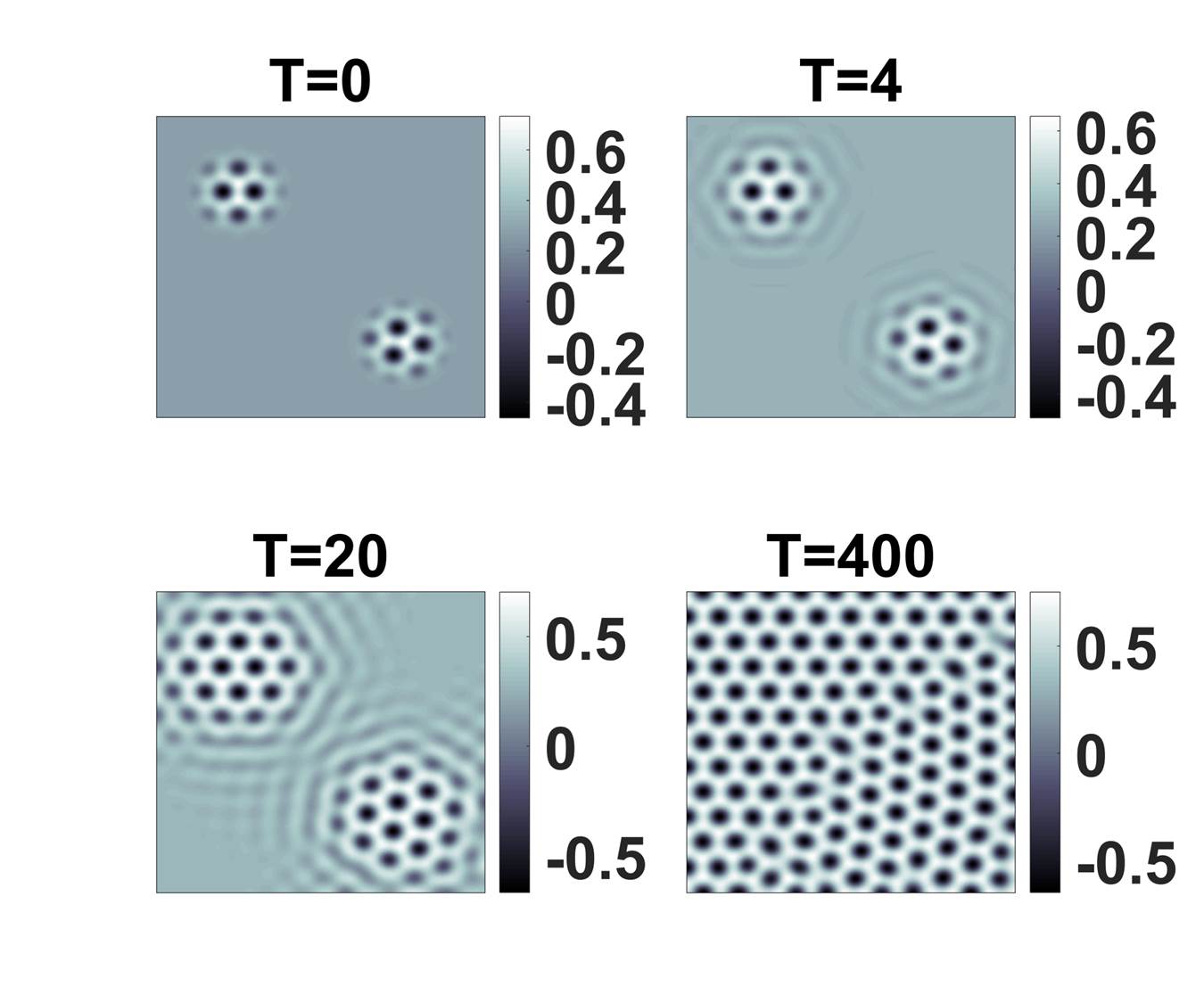}
\end{minipage}}
\subfigure[]{
\begin{minipage}[b]{0.49\linewidth}
\includegraphics[width=1\linewidth]{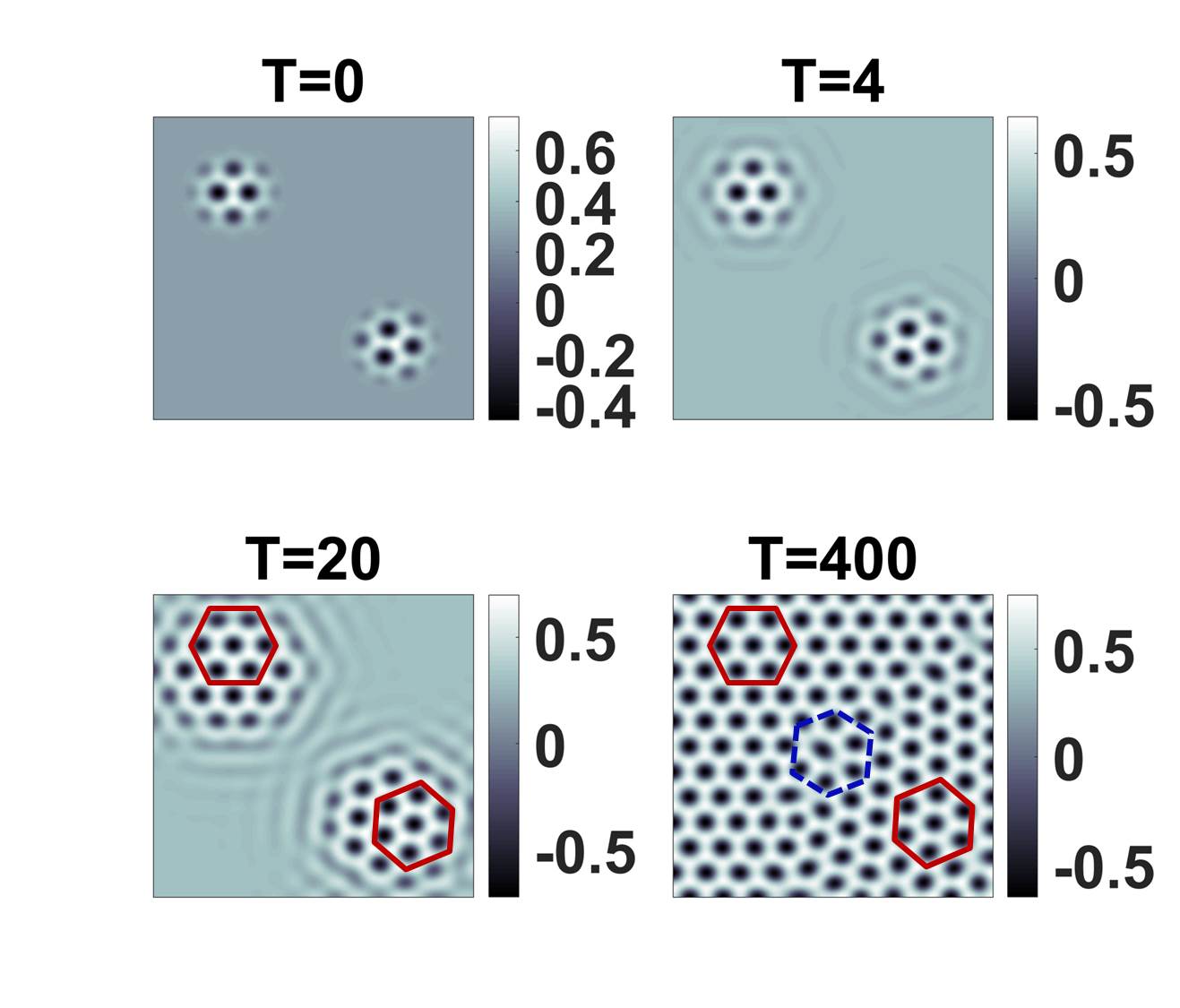}
\end{minipage}}
\subfigure[]{
\begin{minipage}[b]{0.49\linewidth}
\includegraphics[width=1\linewidth]{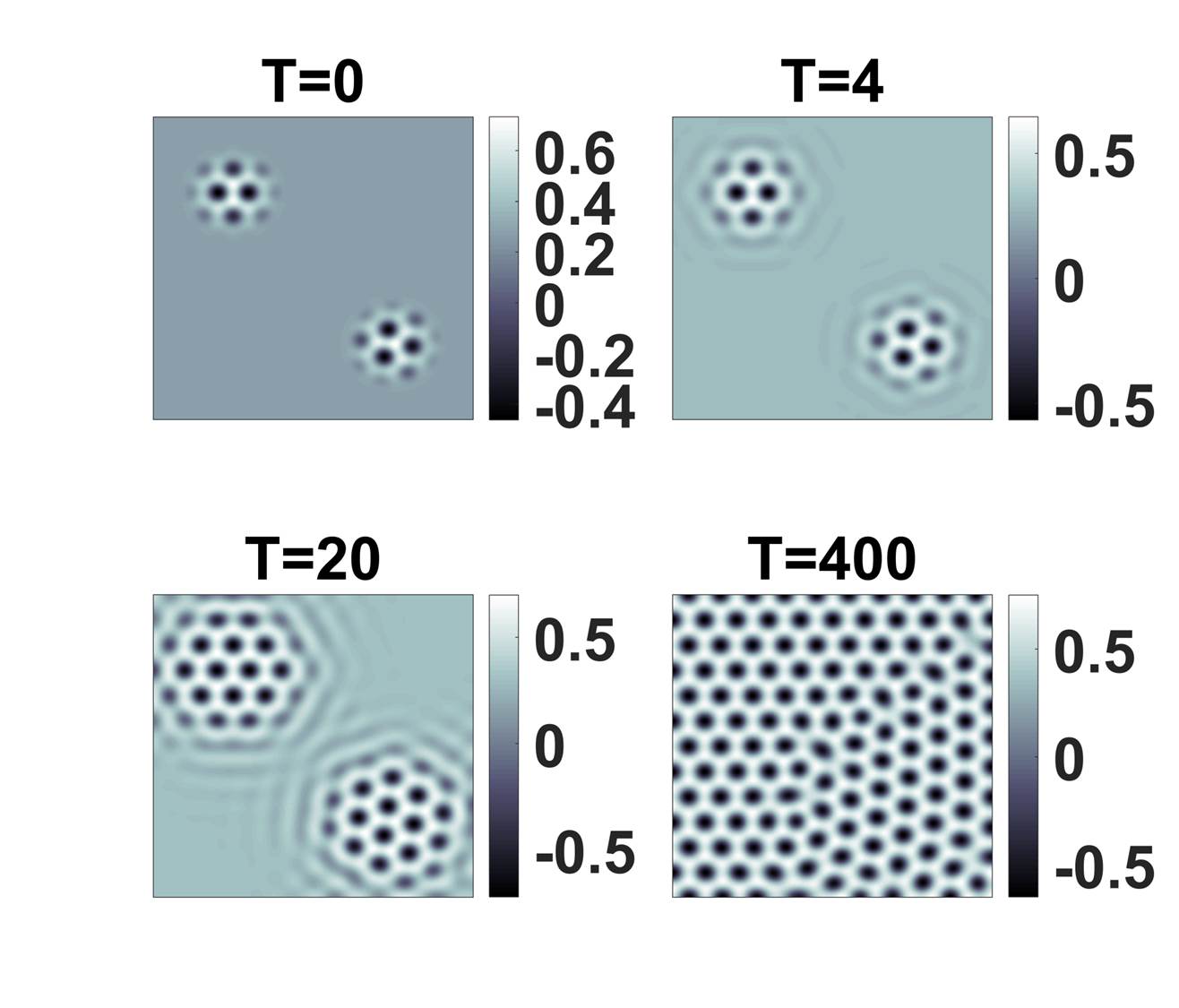}
\end{minipage}}
\caption{Polycrystalline growth with the grain boundary effect simulated using the Allen-Cahn, Cahn-Hillard and Allen-Cahn models with nonlocal constraints at $M=1$.  (a)-(d) are computed using AC-SAV, CH-SAV, AC-P-SAV and AC-L-SAV scheme, respectively. Snapshots of the atomistic density field $\phi$ are taken at $\mathrm T=0, 4, 20, 400$, respectively. Different  growth patterns start affecting each other at $T=20$ and the hexagonal ordering breaks down eventually at the interface of the polyscrystalline. The solid hexagonal ordering and dash non-hexagonal ordering are shown in (c).  Parameter $\eta$ is set as $1\times 10^3$. We use time step $\Delta t=1 \times 10^{-3}$ and  $256\times 256$ space meshes in the 2D simulation. }\label{Fig3}
\end{figure*}

\begin{figure*}
\centering
\subfigure[]{
\begin{minipage}[b]{0.49\linewidth}
\includegraphics[width=1\linewidth]{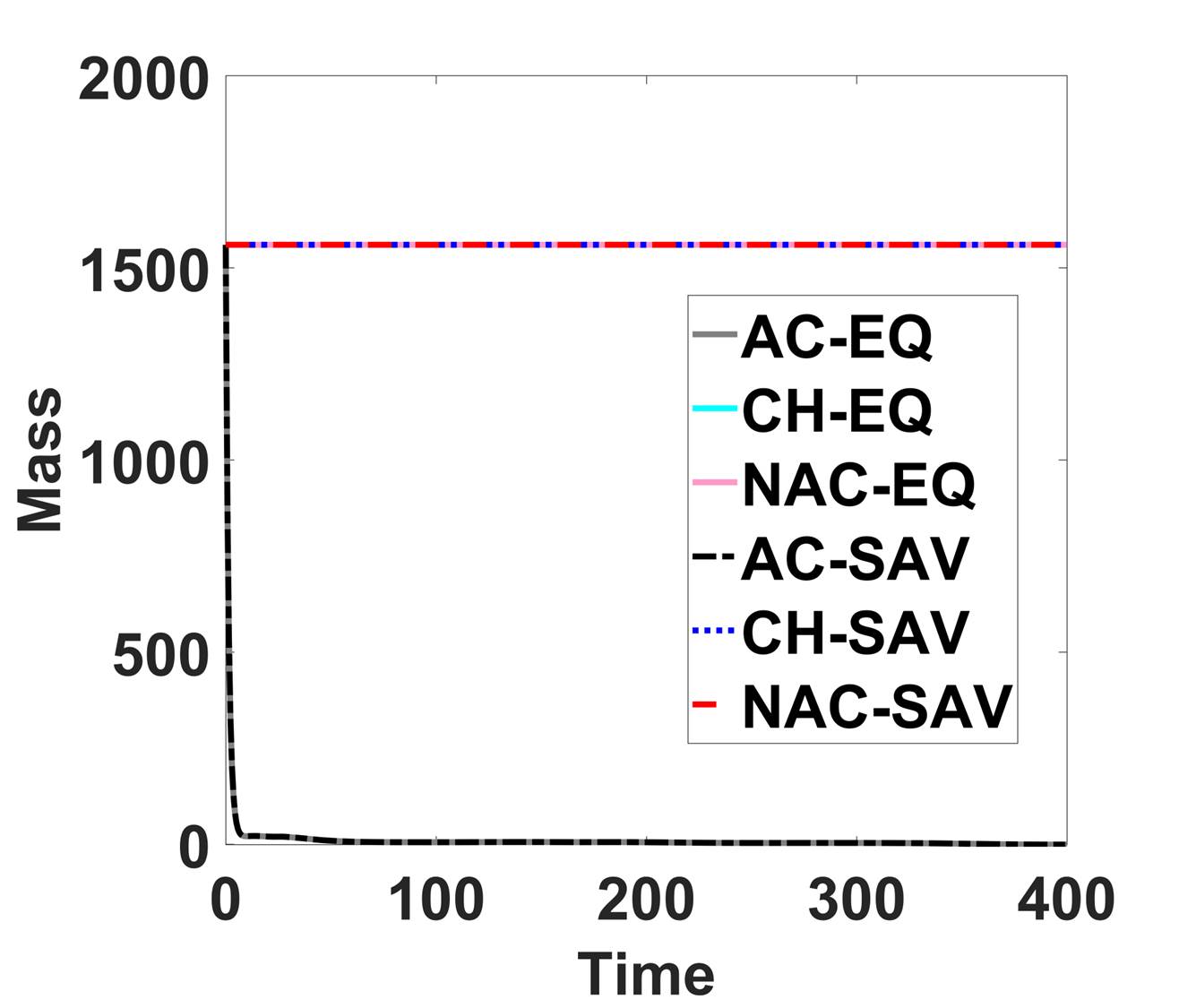}
\end{minipage}}
\subfigure[]{
\begin{minipage}[b]{0.49\linewidth}
\includegraphics[width=1\linewidth]{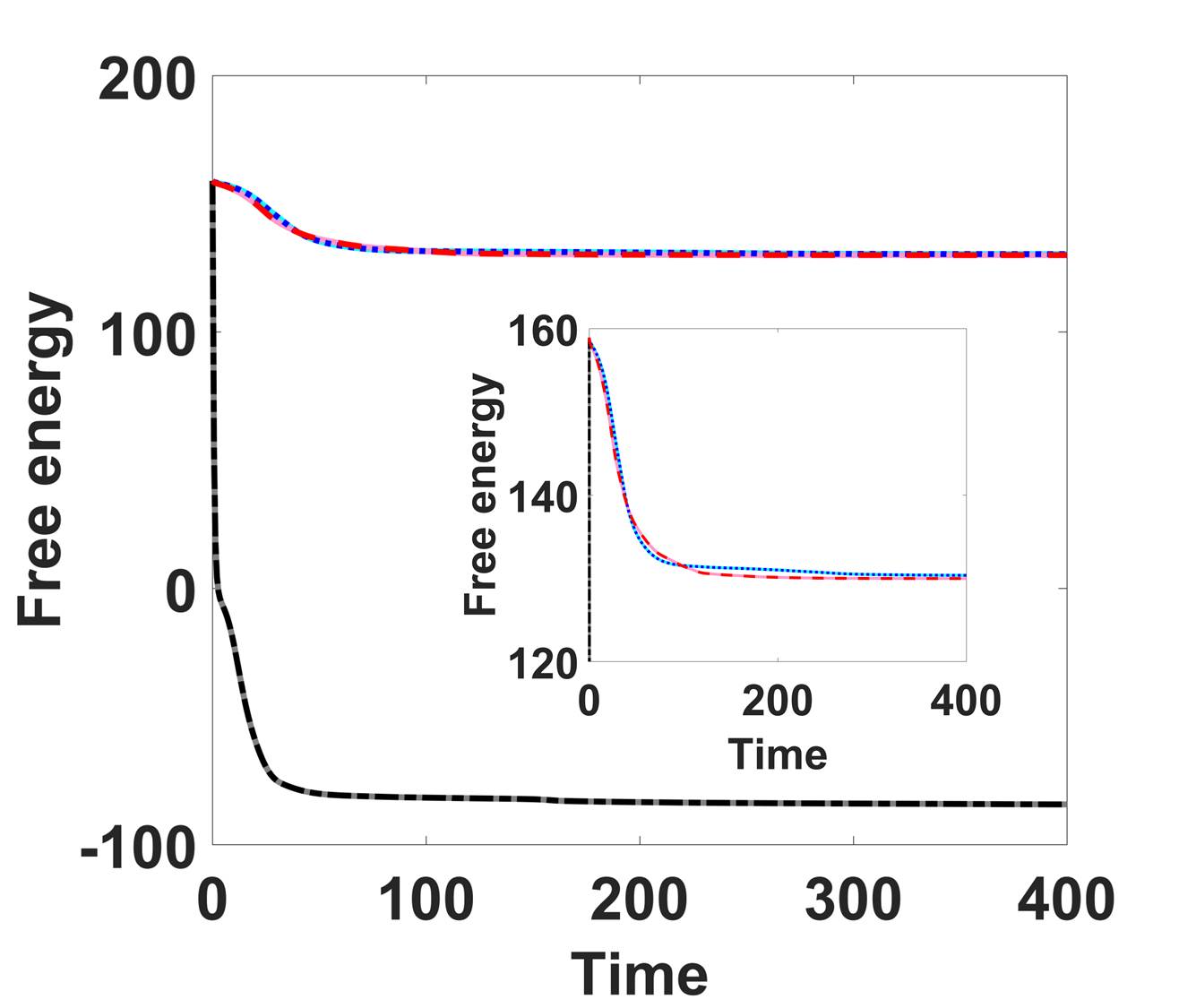}
\end{minipage}}
\caption{Time evolution of mass and free energy from the simulation for Allen-Cahn (AC), Cahn-Hillard (CH) and nonlocal Allen-Cahn (NAC) models are shown in (a) and (b), respectively. Since the dynamical behavior of the Allen-Cahn model with nonlocal constraints are the same, we only show the time evolution of mass and free energy computed using the AC-L-EQ scheme. We compare the time evolution of mass and free energy computed by the EQ and SAV schemes in (a) and (b). We don't see any differences between the results computed by the EQ and SAV methods. The mass computed using the Allen-Cahn model vanishes before $\mathrm T=50$, whereas it is conserved by the Cahn-Hillard and the Allen-Cahn model with nonlocal constraints. Free energy computed by all models are dissipative. Similarly, the nonlocal Allen-Cahn and Cahn-Hillard model predict comparable results in the free energy, but the free energy in the nonlocal models reach the steady state faster than the Cahn-Hillard model does.}\label{Fig4}
\end{figure*}

 Figure \ref{Fig3} depicts the grain boundary effect during polycrystalline growth computed by AC-SAV, CH-SAV, AC-P-SAV and AC-L-SAV schemes, respectively. The snapshots of the phase transitions computed by the Cahn-Hillard and the nonlocal Allen-Cahn models show the Hexagonal ordering are broken at the center at $T=400$. The time evolution of mass and free energy are shown in Figure \ref{Fig4}-(a) and (b), respectively. The mass   in the Allen-Cahn model decays in time. In contrast, the mass in the Cahn-Hillard model and the nonlocal Allen-Cahn models are conserved during the simulation. Meanwhile, the free energies in the Cahn-Hillard model and the nonlocal Allen-Cahn models are larger than that of the  Allen-Cahn model. In Figure \ref{Fig4}-(b), the nonlocal Allen-Cahn models and the Cahn-Hillard model predict comparable time evolution of the free energy but the nonlocal Allen-Cahn models reach the steady state first. The results computed by the EQ and SAV schemes for the same  model (the Allen-Cahn model, Cahn-Hillard model and nonlocal Allen-Cahn model) are nearly identical. We don't see any differences between the results of the Allen-Cahn model with a penalizing potential and the Allen-Cahn model with a Lagrangian multiplier either.

\clearpage

 We note that the results computed by the Allen-Cahn model with a penalizing potential and the Allen-Cahn model with a Lagrangian multiplier in the above two examples at $\eta=1\times 10^3$ are the same. However, the choice of $\eta$ can certainly affect the outcome. In principle, $\eta$ should be chosen as large as possible. However, when $\eta$ is too large, the governing equation becomes very stiff, which forces one to use extremely small time-step size in order to resolve the detail correctly. On the other hand, as we have shown in the two examples, $\eta=1\times 10^3$ is good enough to produce the results that conserve mass very well.

\begin{figure*}
\centering
\subfigure[]{
\begin{minipage}[b]{0.49\linewidth}
\includegraphics[width=1\linewidth]{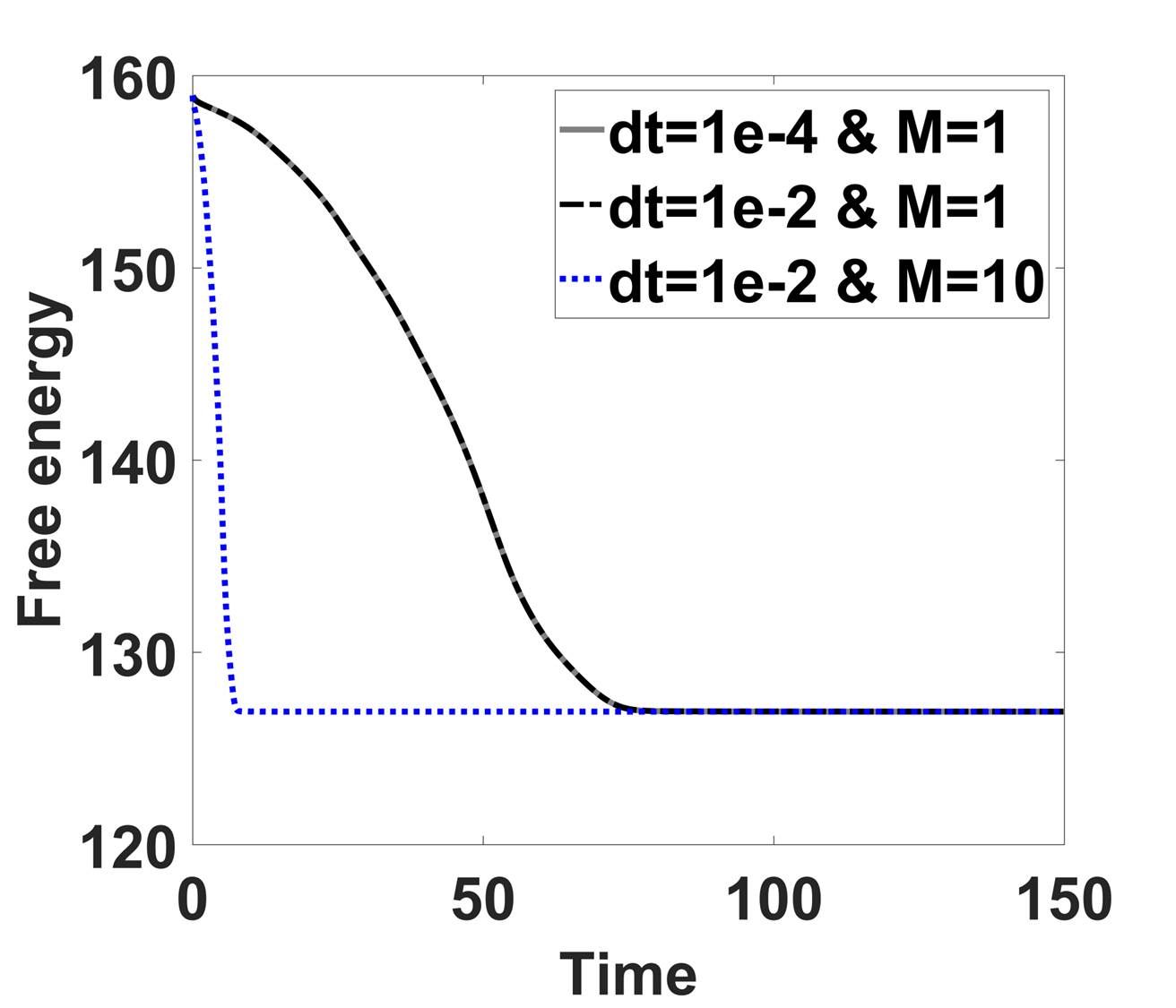}
\end{minipage}}
\subfigure[]{
\begin{minipage}[b]{0.49\linewidth}
\includegraphics[width=1\linewidth]{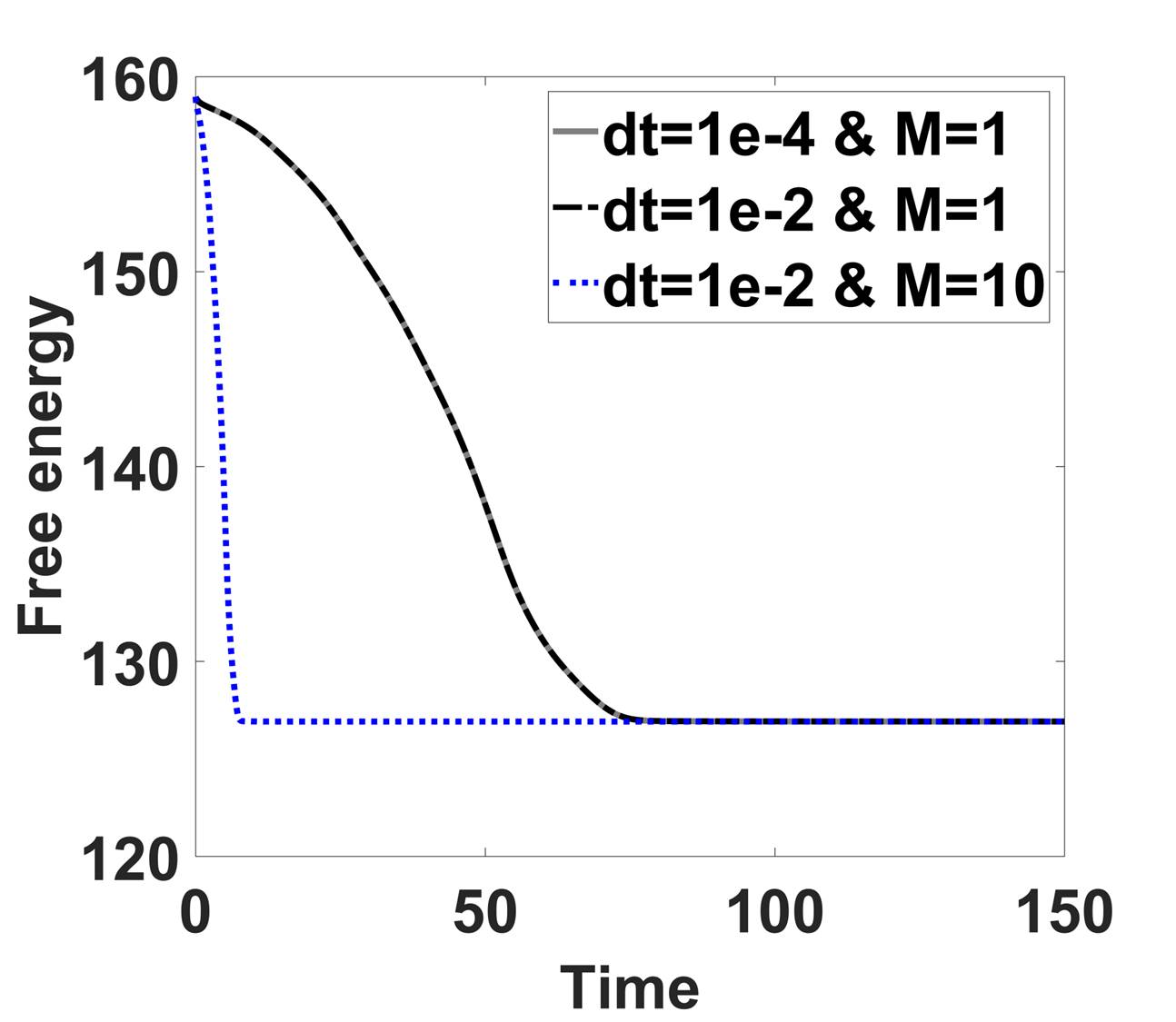}
\end{minipage}}
\subfigure[]{
\begin{minipage}[b]{0.49\linewidth}
\includegraphics[width=1\linewidth]{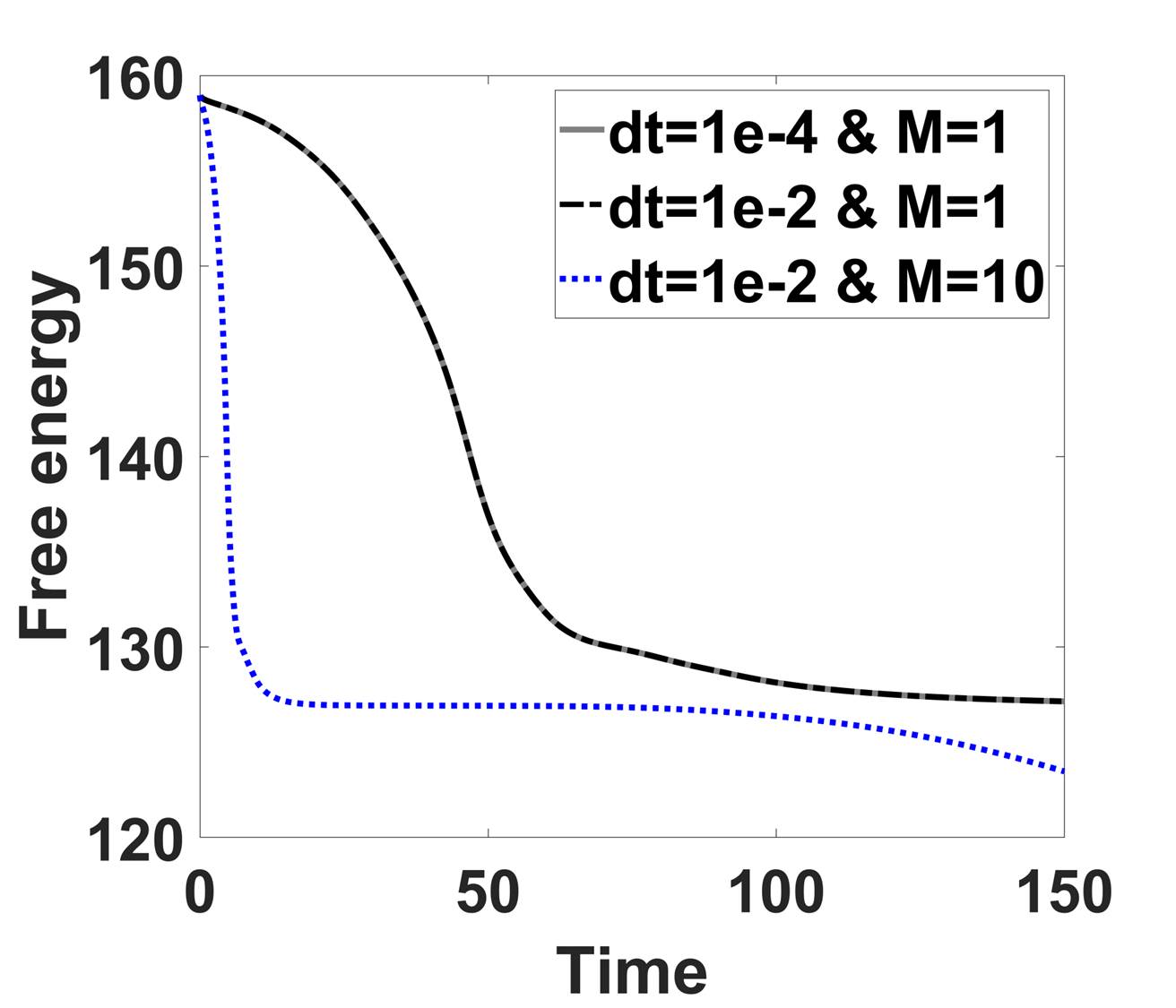}
\end{minipage}}
\subfigure[]{
\begin{minipage}[b]{0.49\linewidth}
\includegraphics[width=1\linewidth]{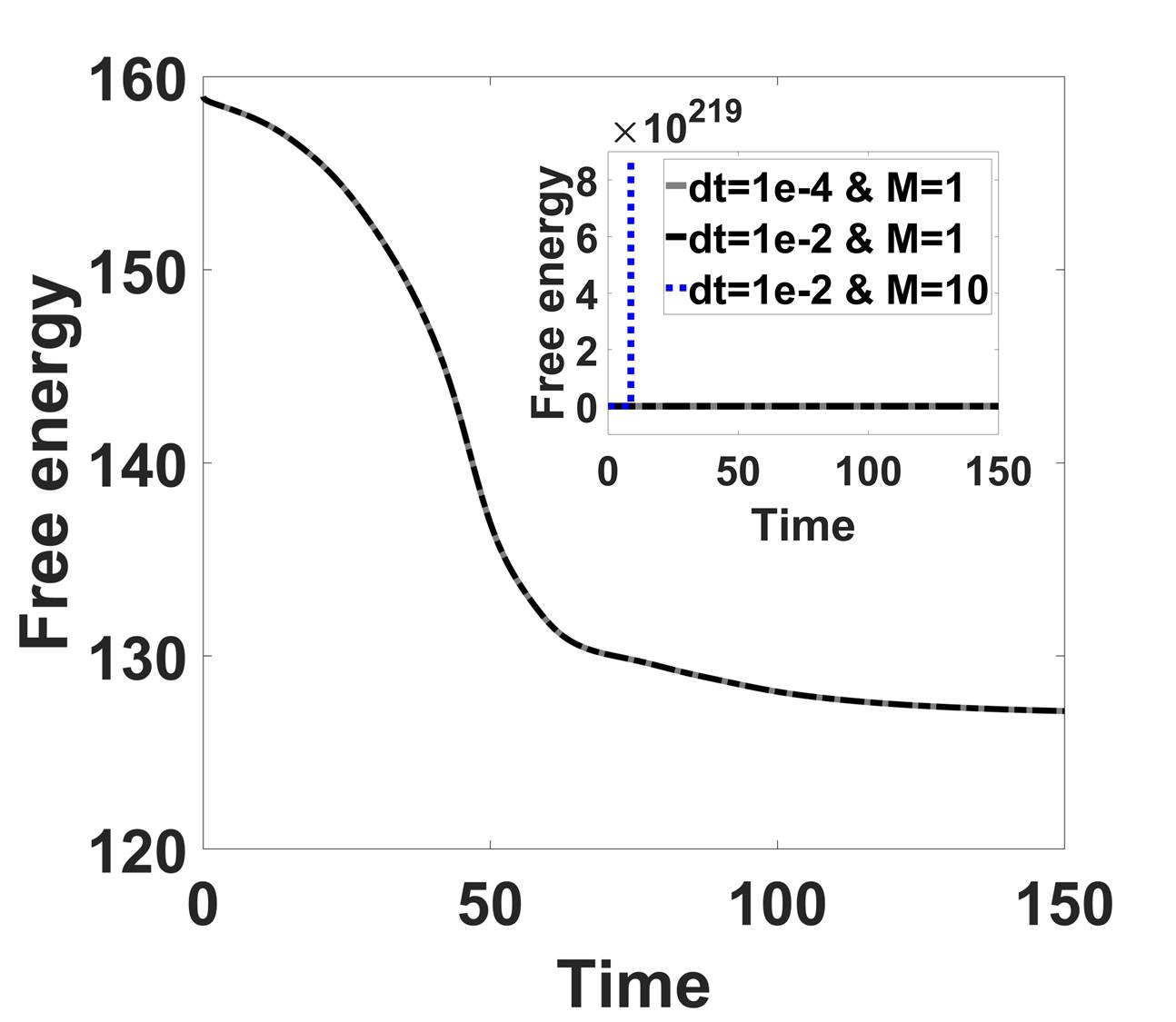}
\end{minipage}}
\caption{Comparison among the nonlocal Allen-Cahn models and the Cahn-Hillard model with large time step size and large mobility coefficients. (a)-(d) are simulated by AC-L-EQ, AC-L-SAV, CH-EQ, CH-SAVscheme with $\Delta t=1\times 10^{-4}$, $M=1$, $\Delta t=1\times 10^{-2}$, $M=1$ and $\Delta t=1\times 10^{-2}$, $M=10$, respectively.  The results show that the nonlocal Allen-Cahn model with a large time step size and large mobility coefficient (a) and (b) performs better than that of the Cahn-Hillard model (c) and (d). The initial conditions and other parameters are chosen the same as those in Figure \ref{Fig1}. }\label{Fig6}
\end{figure*}

Both Figure \ref{Fig2} and Figure \ref{Fig4} show that   the free energy computed by the nonlocal Allen-Cahn models
 reach the steady state faster than that of the Cahn-Hillard model, which is different from  our previous results on a different free energy functional \cite{jing2018second}. Especially, the simulations computed by the  models with mobility coefficient $M=1$ and time steps $1\times 10^{-2},  1\times 10^{-4}$ show the same time evolution behavior. When the mobility is large, say $M=10$, and $\Delta t=0.01$,  the CH-SAV produce an erroneous result while the others produce comparable ones. In practice, if the steady state is more important than the transition dynamics, one can enlarge the mobility coefficient of the   models to accelerate the convergence to  steady state. But, CH-SAV seems to have some accuracy issues with this approach.  The results of the nonlocal Allen-Cahn models in \ref{Fig6}-(a) and (b) show better performance than the Cahn-Hillard model dose in Figure \ref{Fig6}-(c) and (d).

\section{Conclusions}

\noindent \indent We have developed a set of linear, second order, energy stable schemes for the  Allen-Cahn equation with nonlocal constraints that conserve mass and compared them with the energy stable, linear schemes for the Allen-Cahn and the Cahn-Hilliard model. These schemes are devised based on the energy quadratization strategy in the form of EQ and SAV formulation, respectively. We show that they are unconditionally energy stable and uniquely solvable. All schemes can be solved using efficient numerical methods, making the models  alternatives to the Cahn-Hilliard model to describe interface dynamics of immiscible materials while conserving mass. The nonlocal Allen-Cahn models show a faster coarsening rate than that of the Cahn-Hilliard model at the same mobility, but one can enlarge the mobility coefficient of the nonlocal Allen-Cahn model to accelerate their dynamics in case only steady states are of interest. In addition, we have compared the two  Allen-Cahn models with nonlocal constraints numerically. The computational efficiency  of the Allen-Cahn model with a penalizing potential is slightly better than the one with a Lagrange multiplier, but the accuracy of the former depends on a suitable choice of model parameter $\eta$. If the steady state is desired rather than the transient dynamical behavior, large time step size and mobility coefficient can be used to accelerate the simulation. In the end, we show that the  nonlocal Allen-Cahn models perform better that the Cahn-Hilliard model  in the case of a large time step and mobility coefficient.

\section*{Acknowledgements}
 Xiaobo Jing and Qi Wang's research is partially supported by NSFC awards \#11571032, \#91630207 and NSAF-U1530401.
 \section*{Appendix}
\begin{appendices}
\section{Shermann-Morrison formula and its application to solving integro-differential equations}
\noindent \indent Here we give a brief review over the Sherman-Morrison formula \cite{bodewig1959matrix} and explain its applications in the practical implementation of our various relevant schemes.

Suppose $A$ is an invertible square matrix, and $u$,$v$ are column vectors. Then $A+uv^T$ is invertible iff $1+v^TA^{-1}u \neq 0$. If $A+uv^T$ is invertible, then its inverse is given by
\bena
(A+uv^T)^{-1}=A^{-1}-\frac{A^{-1}uv^TA^{-1}}{1+v^TA^{-1}u}.
\eena
So
if $Ay=b$ and $Az=u$, $(A+uv^T)x=b$ has the solution given by
\bena
x=y-\frac{v^Ty}{1+v^Tz}z.
\eena

For the integral term(s) in the  semi-discrete schemes in this study such as (\ref{AC-L-EQ}), we need to discretize it properly. $\forall f$, we discretize $\int_\Omega f \mathrm{d{\bf r}}$ using the composite trapezoidal rule and adding all the elements of the new matrix $w_1 w_2^T f$, where $w_1=\frac{h_x}{2}S$, $w_2=\frac{h_y}{2}S$, $h_x$, $h_y$ are the spatial step sizes and $S={[1,2,2,...,2,2,1]}^T$. For convenience, we use $w_1w_2^T f$ to represent the integral discretized by the composite trapezoidal rule.

To solve equation (\ref{AC-L-EQ}), we  discretize the integral
 or the inner product of functions $(c,\phi^{n+1})d $ as $ u { v}^T \phi^{n+1}$. The scheme is recast to $A\phi^{n+1}+u { v}^T \phi^{n+1}=b^n$.  After using the Sherman-Morrison formula, we get
\bena
\phi ^{n+1}= A^{-1} {b}^n-\frac{{v}^T{ A}^{-1} { {b}}^n}{{1+{ v}^T  A}^{-1}  u}{ A}^{-1}  u,
\eena
In the  inner product of  vectors,
(\ref{AC-L-EQ}) can be rewritten into
\bena
{\phi} ^{n+1}= A^{-1}{b}^n-\frac{\langle c,A^{-1}b^n\rangle }{1+\langle c,A^{-1}d\rangle }{A}^{-1}d.
\eena
So, indeed the approach we take in the study using the discrete inner product is essentially equivalent to applying the Sherman-Morrison formula.

\end{appendices}
\bibliographystyle{plain}
\bibliography{mybibtex}

\begin{thebibliography}{10}

\bibitem{asadi2015review}
Ebrahim Asadi and Mohsen~Asle Zaeem.
\newblock A review of quantitative phase-field crystal modeling of
  solid--liquid structures.
\newblock {\em Jom}, 67(1):186--201, 2015.

\bibitem{backofen2007nucleation}
Rainer Backofen, Andreas R{\"a}tz, and Axel Voigt.
\newblock Nucleation and growth by a phase field crystal (pfc) model.
\newblock {\em Philosophical Magazine Letters}, 87(11):813--820, 2007.

\bibitem{backofen2010phase}
Rainer Backofen and Axel Voigt.
\newblock A phase-field-crystal approach to critical nuclei.
\newblock {\em Journal of Physics: Condensed Matter}, 22(36):364104, 2010.

\bibitem{berry2008simulation}
Joel Berry, KR~Elder, and Martin Grant.
\newblock Simulation of an atomistic dynamic field theory for monatomic
  liquids: Freezing and glass formation.
\newblock {\em Physical Review E}, 77(6):061506, 2008.

\bibitem{bodewig1959matrix}
E~Bodewig.
\newblock Matrix calculus, north, p17, 1959.

\bibitem{chen2018regularized}
Lizhen Chen, Jia Zhao, and Xiaofeng Yang.
\newblock Regularized linear schemes for the molecular beam epitaxy model with
  slope selection.
\newblock {\em Applied Numerical Mathematics}, 128:139--156, 2018.

\bibitem{chen2017uniquely}
Wenbin Chen, Daozhi Han, and Xiaoming Wang.
\newblock Uniquely solvable and energy stable decoupled numerical schemes for
  the cahn--hilliard--stokes--darcy system for two-phase flows in karstic
  geometry.
\newblock {\em Numerische Mathematik}, 137(1):229--255, 2017.

\bibitem{dong2018family}
Suchuan Dong, Zhiguo Yang, and Lianlei Lin.
\newblock A family of second-order energy-stable schemes for cahn-hilliard type
  equations.
\newblock {\em arXiv preprint arXiv:1803.06047}, 2018.

\bibitem{du2018stabilized}
Qiang Du, Lili Ju, Xiao Li, and Zhonghua Qiao.
\newblock Stabilized linear semi-implicit schemes for the nonlocal
  cahn--hilliard equation.
\newblock {\em Journal of Computational Physics}, 363:39--54, 2018.

\bibitem{du2004phase}
Qiang Du, Chun Liu, and Xiaoqiang Wang.
\newblock A phase field approach in the numerical study of the elastic bending
  energy for vesicle membranes.
\newblock {\em Journal of Computational Physics}, 198(2):450--468, 2004.

\bibitem{elder2004modeling}
KR~Elder and Martin Grant.
\newblock Modeling elastic and plastic deformations in nonequilibrium
  processing using phase field crystals.
\newblock {\em Physical Review E}, 70(5):051605, 2004.

\bibitem{elder2002modeling}
KR~Elder, Mark Katakowski, Mikko Haataja, and Martin Grant.
\newblock Modeling elasticity in crystal growth.
\newblock {\em Physical review letters}, 88(24):245701, 2002.

\bibitem{elder2007phase}
KR~Elder, Nikolas Provatas, Joel Berry, Peter Stefanovic, and Martin Grant.
\newblock Phase-field crystal modeling and classical density functional theory
  of freezing.
\newblock {\em Physical Review B}, 75(6):064107, 2007.

\bibitem{elliott1993global}
Charles~M Elliott and AM~Stuart.
\newblock The global dynamics of discrete semilinear parabolic equations.
\newblock {\em SIAM journal on numerical analysis}, 30(6):1622--1663, 1993.

\bibitem{elsey2013simple}
Matt Elsey and Benedikt Wirth.
\newblock A simple and efficient scheme for phase field crystal simulation?
\newblock {\em ESAIM: Mathematical Modelling and Numerical Analysis},
  47(5):1413--1432, 2013.

\bibitem{eyre1998unconditionally}
David~J Eyre.
\newblock Unconditionally gradient stable time marching the cahn-hilliard
  equation.
\newblock {\em MRS Online Proceedings Library Archive}, 529, 1998.

\bibitem{fan2017componentwise}
Xiaolin Fan, Jisheng Kou, Zhonghua Qiao, and Shuyu Sun.
\newblock A componentwise convex splitting scheme for diffuse interface models
  with van der waals and peng--robinson equations of state.
\newblock {\em SIAM Journal on Scientific Computing}, 39(1):B1--B28, 2017.

\bibitem{gomez2011provably}
Hector Gomez and Thomas~JR Hughes.
\newblock Provably unconditionally stable, second-order time-accurate, mixed
  variational methods for phase-field models.
\newblock {\em Journal of Computational Physics}, 230(13):5310--5327, 2011.

\bibitem{gomez2012unconditionally}
Hector Gomez and Xes{\'u}s Nogueira.
\newblock An unconditionally energy-stable method for the phase field crystal
  equation.
\newblock {\em Computer Methods in Applied Mechanics and Engineering},
  249:52--61, 2012.

\bibitem{Gong2018Linear}
Yuezheng Gong, Jia Zhao, and Qi~Wang.
\newblock Linear second order in time energy stable schemes for hydrodynamic
  models of binary mixtures based on a spatially pseudospectral approximation.
\newblock {\em Advances in Computational Mathematics}, pages 1--28, 2018.

\bibitem{Gong2018Fully}
Yuezheng Gong, Jia Zhao, Xiaogang Yang, and Qi~Wang.
\newblock Fully discrete second-order linear schemes for hydrodynamic phase
  field models of binary viscous fluid flows with variable densities.
\newblock {\em Siam Journal on Scientific Computing}, 40(1):B138--B167, 2018.

\bibitem{guttenberg2010emergence}
Nicholas Guttenberg, Nigel Goldenfeld, and Jonathan Dantzig.
\newblock Emergence of foams from the breakdown of the phase field crystal
  model.
\newblock {\em Physical Review E}, 81(6):065301, 2010.

\bibitem{han2015second}
Daozhi Han and Xiaoming Wang.
\newblock A second order in time, uniquely solvable, unconditionally stable
  numerical scheme for cahn--hilliard--navier--stokes equation.
\newblock {\em Journal of Computational Physics}, 290:139--156, 2015.

\bibitem{hu2009stable}
Zhengzheng Hu, Steven~M Wise, Cheng Wang, and John~S Lowengrub.
\newblock Stable and efficient finite-difference nonlinear-multigrid schemes
  for the phase field crystal equation.
\newblock {\em Journal of Computational Physics}, 228(15):5323--5339, 2009.

\bibitem{jaatinen2010extended}
A~Jaatinen and T~Ala-Nissila.
\newblock Extended phase diagram of the three-dimensional phase field crystal
  model.
\newblock {\em Journal of Physics: Condensed Matter}, 22(20):205402, 2010.

\bibitem{jing2018second}
Xiaobo Jing, Jun Li, Xueping Zhao, and Qi~Wang.
\newblock Second order linear energy stable schemes for allen-cahn equations
  with nonlocal constraints.
\newblock {\em arXiv preprint arXiv:1810.05311}, 2018.

\bibitem{karma1996phase}
Alain Karma and Wouter-Jan Rappel.
\newblock Phase-field method for computationally efficient modeling of
  solidification with arbitrary interface kinetics.
\newblock {\em Physical review E}, 53(4):R3017, 1996.

\bibitem{li2018unconditionally}
Hongwei Li, Lili Ju, Chenfei Zhang, and Qiujin Peng.
\newblock Unconditionally energy stable linear schemes for the diffuse
  interface model with peng--robinson equation of state.
\newblock {\em Journal of Scientific Computing}, 75(2):993--1015, 2018.

\bibitem{lowen2010phase}
Hartmut L{\"o}wen.
\newblock A phase-field-crystal model for liquid crystals.
\newblock {\em Journal of Physics: Condensed Matter}, 22(36):364105, 2010.

\bibitem{mellenthin2008phase}
Jesper Mellenthin, Alain Karma, and Mathis Plapp.
\newblock Phase-field crystal study of grain-boundary premelting.
\newblock {\em Physical Review B}, 78(18):184110, 2008.

\bibitem{onsager1931reciprocal1}
Lars Onsager.
\newblock Reciprocal relations in irreversible processes. i.
\newblock {\em Physical review}, 37(4):405, 1931.

\bibitem{onsager1931reciprocal2}
Lars Onsager.
\newblock Reciprocal relations in irreversible processes. ii.
\newblock {\em Physical review}, 38(12):2265, 1931.

\bibitem{pisutha2013calculations}
N~Pisutha-Arnond, VWL Chan, KR~Elder, and K~Thornton.
\newblock Calculations of isothermal elastic constants in the phase-field
  crystal model.
\newblock {\em Physical Review B}, 87(1):014103, 2013.

\bibitem{provatas2007using}
N~Provatas, JA~Dantzig, B~Athreya, P~Chan, P~Stefanovic, N~Goldenfeld, and
  KR~Elder.
\newblock Using the phase-field crystal method in the multi-scale modeling of
  microstructure evolution.
\newblock {\em Jom}, 59(7):83--90, 2007.

\bibitem{provatas2011phase}
Nikolas Provatas and Ken Elder.
\newblock {\em Phase-field methods in materials science and engineering}.
\newblock John Wiley \& Sons, 2011.

\bibitem{Rubinstein1992Nonlocal}
Jacob Rubinstein and Peter Sternberg.
\newblock Nonlocal reaction¡ªdiffusion equations and nucleation.
\newblock {\em Ima Journal of Applied Mathematics}, 48(3):249--264, 1992.

\bibitem{shen2012second}
Jie Shen, Cheng Wang, Xiaoming Wang, and Steven~M Wise.
\newblock Second-order convex splitting schemes for gradient flows with
  ehrlich--schwoebel type energy: application to thin film epitaxy.
\newblock {\em SIAM Journal on Numerical Analysis}, 50(1):105--125, 2012.

\bibitem{shen2018convergence}
Jie Shen and Jie Xu.
\newblock Convergence and error analysis for the scalar auxiliary variable
  (sav) schemes to gradient flows.
\newblock {\em SIAM Journal on Numerical Analysis}, 56(5):2895--2912, 2018.

\bibitem{shen2018scalar}
Jie Shen, Jie Xu, and Jiang Yang.
\newblock The scalar auxiliary variable (sav) approach for gradient flows.
\newblock {\em Journal of Computational Physics}, 353:407--416, 2018.

\bibitem{shen2010numerical}
Jie Shen and Xiaofeng Yang.
\newblock Numerical approximations of allen-cahn and cahn-hilliard equations.
\newblock {\em Discrete Contin. Dyn. Syst}, 28(4):1669--1691, 2010.

\bibitem{tegze2009advanced}
Gy{\"o}rgy Tegze, Gurvinder Bansel, Gyula~I T{\'o}th, Tam{\'a}s Pusztai,
  Zhongyun Fan, and L{\'a}szl{\'o} Gr{\'a}n{\'a}sy.
\newblock Advanced operator splitting-based semi-implicit spectral method to
  solve the binary phase-field crystal equations with variable coefficients.
\newblock {\em Journal of Computational Physics}, 228(5):1612--1623, 2009.

\bibitem{toth2010polymorphism}
Gyula~I T{\'o}th, Gy{\"o}rgy Tegze, Tam{\'a}s Pusztai, Gergely T{\'o}th, and
  L{\'a}szl{\'o} Gr{\'a}n{\'a}sy.
\newblock Polymorphism, crystal nucleation and growth in the phase-field
  crystal model in 2d and 3d.
\newblock {\em Journal of Physics: Condensed Matter}, 22(36):364101, 2010.

\bibitem{vignal2015energy}
Philippe Vignal, Lisandro Dalcin, Donald~L Brown, Nathan Collier, and Victor~M
  Calo.
\newblock An energy-stable convex splitting for the phase-field crystal
  equation.
\newblock {\em Computers \& Structures}, 158:355--368, 2015.

\bibitem{wang2011energy}
C~Wang and Steven~M Wise.
\newblock An energy stable and convergent finite-difference scheme for the
  modified phase field crystal equation.
\newblock {\em SIAM Journal on Numerical Analysis}, 49(3):945--969, 2011.

\bibitem{wang2010unconditionally}
Cheng Wang, Xiaoming Wang, and Steven~M Wise.
\newblock Unconditionally stable schemes for equations of thin film epitaxy.
\newblock {\em Discrete Contin. Dyn. Syst}, 28(1):405--423, 2010.

\bibitem{wang2017efficient}
Lin Wang and Haijun Yu.
\newblock On efficient second order stabilized semi-implicit schemes for the
  cahn--hilliard phase-field equation.
\newblock {\em Journal of Scientific Computing}, pages 1--25, 2017.

\bibitem{wise2009energy}
Steven~M Wise, Cheng Wang, and John~S Lowengrub.
\newblock An energy-stable and convergent finite-difference scheme for the
  phase field crystal equation.
\newblock {\em SIAM Journal on Numerical Analysis}, 47(3):2269--2288, 2009.

\bibitem{yang2018efficient}
Xiaofeng Yang.
\newblock Efficient schemes with unconditionally energy stability for the
  anisotropic cahn-hilliard equation using the stabilized-scalar augmented
  variable (s-sav) approach.
\newblock {\em arXiv preprint arXiv:1804.02619}, 2018.

\bibitem{yang2006numerical}
Xiaofeng Yang, James~J Feng, Chun Liu, and Jie Shen.
\newblock Numerical simulations of jet pinching-off and drop formation using an
  energetic variational phase-field method.
\newblock {\em Journal of Computational Physics}, 218(1):417--428, 2006.

\bibitem{yang2017linearly}
Xiaofeng Yang and Daozhi Han.
\newblock Linearly first-and second-order, unconditionally energy stable
  schemes for the phase field crystal model.
\newblock {\em Journal of Computational Physics}, 330:1116--1134, 2017.

\bibitem{Yang2017Linear}
Xiaofeng Yang and Lili Ju.
\newblock Linear and unconditionally energy stable schemes for the binary
  fluid¨csurfactant phase field model.
\newblock {\em Computer Methods in Applied Mechanics and Engineering},
  318:1005--1029, 2017.

\bibitem{yang2017numerical}
Xiaofeng Yang, Jia Zhao, and Qi~Wang.
\newblock Numerical approximations for the molecular beam epitaxial growth
  model based on the invariant energy quadratization method.
\newblock {\em Journal of Computational Physics}, 333:104--127, 2017.

\bibitem{Zhaoetal2018-2}
Jia Zhao, Yuezheng Gong, and Qi~Wang.
\newblock Aritrary high order unconditionally energy stable schemes for
  gradient flow models.
\newblock {\em Journal of Computational Physics}, 2018.

\bibitem{Zhao2016Numerical}
Jia Zhao, Qi~Wang, and Xiaofeng Yang.
\newblock Numerical approximations to a new phase field model for two phase
  flows of complex fluids.
\newblock {\em Computer Methods in Applied Mechanics and Engineering},
  310:77--97, 2016.

\bibitem{Zhao2017Numerical}
Jia Zhao, Qi~Wang, and Xiaofeng Yang.
\newblock Numerical approximations for a phase field dendritic crystal growth
  model based on the invariant energy quadratization approach.
\newblock {\em International Journal for Numerical Methods in Engineering},
  110(3), 2017.

\bibitem{zhao2017novel}
Jia Zhao, Xiaofeng Yang, Yuezheng Gong, and Qi~Wang.
\newblock A novel linear second order unconditionally energy stable scheme for
  a hydrodynamic-tensor model of liquid crystals.
\newblock {\em Computer Methods in Applied Mechanics and Engineering},
  318:803--825, 2017.

\bibitem{zhao2018general}
Jia Zhao, Xiaofeng Yang, Yuezheng Gong, Xueping Zhao, Xiaogang Yang, Jun Li,
  and Qi~Wang.
\newblock A general strategy for numerical approximations of non-equilibrium
  models--part i: Thermodynamical systems.
\newblock {\em International Journal of Numerical Analysis \& Modeling},
  15(6):884--918, 2018.

\bibitem{zhao2016energy}
Jia Zhao, Xiaofeng Yang, Jun Li, and Qi~Wang.
\newblock Energy stable numerical schemes for a hydrodynamic model of nematic
  liquid crystals.
\newblock {\em SIAM Journal on Scientific Computing}, 38(5):A3264--A3290, 2016.

\bibitem{zhao2016decoupled}
Jia Zhao, Xiaofeng Yang, Jie Shen, and Qi~Wang.
\newblock A decoupled energy stable scheme for a hydrodynamic phase-field model
  of mixtures of nematic liquid crystals and viscous fluids.
\newblock {\em Journal of Computational Physics}, 305:539--556, 2016.

\end{thebibliography}
\end{document}